\documentclass[
  a4paper,
  11pt,
  oneside
]{article}

\usepackage[utf8]{inputenc}
\usepackage{lmodern}
\usepackage[german,english]{babel}

\usepackage[numbers]{natbib}
\usepackage[linesnumbered,ruled,vlined]{algorithm2e}
\usepackage{listings}

\lstdefinestyle{cppstyle}{
  backgroundcolor   = \color{white},
  commentstyle      = \color{gray},
  keywordstyle      = \color{LUH-blue},
  numberstyle       = \tiny\color{codegray},
  stringstyle       = \color{codepurple},
  basicstyle        = \ttfamily\footnotesize,
  breakatwhitespace = false,
  breaklines        = true,
  captionpos        = b,
  frame             = leftline,
  keepspaces        = true,
  numbers           = left,
  numbersep         = 5pt,
  showspaces        = false,
  showstringspaces  = false,
  showtabs          = false,
  tabsize           = 2
}

\usepackage{enumerate}

\usepackage{graphicx}   % needed for graphics and tikz
\usepackage{float}      % figure enviroment
\usepackage{wrapfig}    % figures wrapped with text
\usepackage{caption}    % captions for figures
\usepackage{tikz}       % Tikz ist kein Zeichen Programm
  \usetikzlibrary{arrows.meta}
  \usetikzlibrary{calc}

\usepackage{pgf}
\usepackage{pgfplots}

\usepackage{multirow}
\usepackage{rotating}

\usepackage{booktabs}  % toprule, midrule, buttomrule

\usepackage{url}
\usepackage{authblk}
\usepackage{hyperref}
  \hypersetup{
    colorlinks,
    citecolor=black,
    filecolor=black,
    linkcolor=black,
    urlcolor=black
  }

\usepackage{mathtools}  % implements: \coloneqq :=
\usepackage{xfrac}      % implements: \sfrac 1/2
\usepackage{amssymb}    % to use the mathbb enviroment
\usepackage{amsmath}    
\usepackage{amsfonts}
\usepackage{units}
\numberwithin{equation}{section}
\numberwithin{figure}{section}
\numberwithin{table}{section}
\usepackage{amscd}

\usepackage{microtype}  % better text setting
\usepackage{relsize}    % resize text
\usepackage[            % margins
  top=2.5cm, 
  bottom=2.0cm, 
  left=2.5cm, 
  right=2.5cm
]{geometry}
\usepackage{xcolor}
\definecolor{LUH-blue} 	  {RGB}{ 38,  87, 167}
\definecolor{LUH-lblue} 	{RGB}{ 92, 129, 189}
\definecolor{LUH-llblue}	{RGB}{146, 170, 210}
\definecolor{LUH-lllblue}	{RGB}{204, 220, 235}

\definecolor{PXD-grey}		{RGB}{217, 217, 217}
\definecolor{PXD-lgrey}		{RGB}{227, 227, 227}
\definecolor{PXD-llgrey}	{RGB}{235, 235, 235}

\definecolor{PXD-yellow}	{RGB}{255, 205,   9}
\definecolor{PXD-lyellow}	{RGB}{255, 218,  70}
\definecolor{PXD-llyellow}{RGB}{255, 229, 131}

\definecolor{PXD-red}		  {RGB}{197,  18,  88}
\definecolor{PXD-lred}		{RGB}{212,  77, 130}
\definecolor{PXD-llred}		{RGB}{225, 136, 171}

\definecolor{LUH-green} 	{RGB}{200, 211,  23}
\definecolor{LUH-lgreen} 	{RGB}{214, 222,  81}
\definecolor{LUH-llgreen} {RGB}{227, 232, 138}

\definecolor{LUH-gray}		{RGB}{153, 153, 153}
\definecolor{LUH-lgray}		{RGB}{204, 204, 204}
\definecolor{LUH-llgray}	{RGB}{235, 235, 235}

\definecolor{dark-green} 	{RGB}{0, 80,  5}

\usepackage{amsthm}
\usepackage{color}

\theoremstyle{plain}

\theoremstyle{definition}

\theoremstyle{remark}

\newcommand{\dealii}{\texttt{deal.II}}

\renewcommand{\vec}[1]{\mathbf{#1}}

\usepackage[linesnumbered,ruled,vlined]{algorithm2e}
\SetKwBlock{Loop}{Loop}{end}
\SetKwBlock{LoopOverAll}{Loop over all}{end}

\begin{document}

%+-------------------------------------------------------------------+
%|                         Headline and Affiliations                 |
%+-------------------------------------------------------------------+
\title{%
  Algorithmic realization of the solution to the sign conflict problem
  for hanging nodes on hp-hexahedral Nédélec elements
}
\author[1,2]{S. Kinnewig}
\author[1,2]{T. Wick}
\author[1,2]{S. Beuchler}

\affil[1]{%
  Leibniz University Hannover,
  Institute of Applied Mathematics,
  Welfengarten 1,
  30167 Hannover,
  Germany
}
\affil[2]{%
  Cluster of Excellence PhoenixD (Photonics, Optics, and
  Engineering - Innovation Across Disciplines),
  Leibniz University Hannover,
  Germany
}

\date{}

\maketitle
	
%+-------------------------------------------------------------------+
%|                            Abstract                               |
%+-------------------------------------------------------------------+
\begin{abstract}
  In this work, Nédélec elements on locally refined meshes with hanging 
  nodes are considered. A crucial aspect is the orientation of the hanging 
  edges and faces. 
  For non-orientable meshes, no solution or implementation has been 
  available to date. The problem statement and corresponding algorithms 
  are described in great detail. As a model problem, the time-harmonic 
  Maxwell's equations are adopted because Nédélec elements constitute 
  their natural discretization. The algorithms and implementation are 
  demonstrated through two numerical examples on different uniformly and 
  adaptively refined meshes. The implementation is performed within the 
  finite element library \dealii.
\end{abstract}
	
%+-------------------------------------------------------------------+
%|                     Section: Introduction                         |
%+-------------------------------------------------------------------+
\section{Introduction}
\label{introduction}
This work is devoted to the numerical construction of the N\'ed\'elec elements 
in three spatial dimensions on locally refined meshes with hanging nodes. 
N\'ed\'elec elements are usually required for solving Maxwell's equations 
\cite{Bk:Fey:63,Bk:Monk:2003}, which are fundamental to many fields of research. 
They have numerous practical applications, ranging from Magnetic Induction Tomography (MIT) in
medicine \cite{Art:Zolgharni:MIT:2009}, geo-electromagnetic modeling in geophysics
\cite{Art:Grayver:Geoelectromagnetic:2015}, to quantum computing \cite{Art:Kues:23},
quantum communication in optics \cite{Art:Melchert:Soliton:2023}, and photonics 
as they are of interest in the cluster of excellence    
PhoenixD\footnote{\url{https://www.phoenixd.uni-hannover.de/en/}}.
Designing optical components can be challenging, and simulations are often necessary for support.
These simulations involve modeling electromagnetic waves within the components, which is achieved
by solving Maxwell's equations using Nédélec elements as the natural finite element (FE) discretization. The discretized systems result in linear equation systems. 
Besides efficient numerical solution schemes (e.g., \cite{Art:Bonazetal:19,Art:Gander:OptimizedSchwarz:12,Art:peccc:2023,Art:FauParMel:23,Art:ErnstGander:12,Art:Arnold:00,Art:Hiptmair:99,Art:Toselli:06,Art:Dohrmann:16}),
solving Maxwell's equations remains computationally expensive.
Therefore, adaptive strategies such as local grid refinement are highly desirable.
These strategies can keep computational costs reasonable while increasing accuracy. 
They can be achieved with heuristic error indicators, geometry-oriented refinement, residual-based error control, or goal-oriented error control. 
The discussion of error estimators is outside the scope of this work, but we refer the reader to
\cite{Art:BaRhei:78,Art:ZiZhu:92b,Bk:Verfuerth:1996,Bk:AinsworthOden:00,Bk:Bangerth:03,Bk:Repin:08,Bk:ErEstHaJoh:09}.

The key objective of this work is to address a long-standing open problem that concerns the design of algorithms and corresponding implementations of the N\'ed\'elec basis functions in three dimensions on non-orientable locally refined meshes.
As previously mentioned, the authors of \cite{Art:Kynch:ResolvingTheSignConflict:2017}
considered high-polynomial N\'ed\'elec basis functions to capture skin effects that appear
in the MIT problem. 
Therefore, they described a procedure to overcome the sign
conflict on $hp$-N\'ed\'elec elements. In \dealii, prior work already utilized
hanging nodes for N\'ed\'elec elements, such as the work of Bürg \cite{Art:Buerg:12}.
However, an older implementation, the so-called \texttt{FE\_Nedelec}\footnote{\url{https://www.dealii.org/current/doxygen/deal.II/classFE__Nedelec.html}} 
was used, and it can only be applied to oriented grids.

Our choice for a suitable programming platform is motivated
by modern available FEM libraries that include support for high-order
N\'ed\'elec elements. Various open-source finite element libraries allow
the use of Nédélec elements of polynomial degree $p\geq2$. The \texttt{Elmer FEM}
library \cite{Software:ELMER} can handle unstructured grids with a maximum of $p=2$,
while \texttt{FreeFEM++} \cite{Software:FreeFEM} can support a maximum of $p=3$.
\texttt{NGSolve} \cite{Software:NGSolve,Software:Netgen} utilized the basis functions introduced by 
Zaglmayr \cite{Diss:Zaglmayr:2006} to implement high polynomial functions on unstructured grids. \texttt{hp3D} \cite{Software:hp3D} implements the Nédélec
functions based on the hierarchical polynomial basis from Demkowicz \cite{Bk:Demkowicz.I:2008}.
Also, the libraries  \texttt{FEniCS}
\cite{Software:FEniCS} (unstructured), \texttt{MFEM} \cite{Software:MFEM}, and \texttt{GetDP}
\cite{Software:GetDP} (unstructured) implement high polynomial Nédélec elements. Moreover, \texttt{GetDDM} \cite{Software:GetDDM} 
is an extension of \texttt{GetDP} that implements optimized Schwarz domain decomposition methods, which is a
well-established method for solving the ill-posed Maxwell's problems.

We have chosen \dealii \cite{Software:dealii:9.5,Art:dealii:20} as it offers high-polynomial (i.e., arbitrary polynomial degrees $p$) N\'ed\'elec basis functions based on Schöberl and Zaglmayr's basis function set for the complete De-Rham sequence \cite{Diss:Zaglmayr:2006}; see also \cite{Art:Ainsworth:01} for the two-dimensional case.
\dealii is well-established, with a large user base and excellent accessibility,
thanks to its comprehensive documentation, which is essential for sustainable software development.
It uses tensor product elements and is designed with adaptive mesh refinement in mind, providing a range
of functionalities for the computation of error estimators.
Due to the use of quadrilateral and hexahedral elements, local mesh refinement in \dealii
requires the use of hanging nodes.
As a starting point for our implementation of hanging nodes, we use the work of Ledger and Kynch
\cite{Art:Kynch:ResolvingTheSignConflict:2017} for non-orientable grids.

In more detail, we extend \dealii's class \texttt{FE\_NedelecSZ}\footnote{\url{https://www.dealii.org/current/doxygen/deal.II/classFE__NedelecSZ.html}}, which can also be applied to non-orientable grids. 
The extension to three dimensions is non-trivial, as we shall see.
The main work here relies upon the high number of possible configurations we have to cope with.
To overcome the sign conflict
in the case of hanging edges and faces, we need to adapt the associated constraint matrix
that restricts the additional Degrees of Freedom (DoFs) introduced by the hanging edges and faces accordingly. 
One face has $2^3$ possible orientations, which results after local refinement 
into four child faces. Consequently, we have to deal with $2^{15}$ possible configurations. 
As dealing with every case individually would be even more cumbersome, 
we treat the outer edges, the inner edges, and the faces separately to reduce the number of necessary algorithms in order to obtain
an efficient code.
Our goal is to resolve sign conflicts regardless of the
polynomial degree involved. To achieve this, we need to comprehend the structure of the
constraint matrix so that we can develop algorithms that can deal with any given polynomial degree.
As one of our aims is to make these results accessible, we provide the most crucial steps
as pseudo-code. 
Our implementation is available open-source at \cite{Software:hanging-nodes-for-nedelec}\footnote{\url{https://zenodo.org/records/10913219}}.
These accomplishments are exemplarily applied to the time-harmonic Maxwell's equations,
which are solved for two different configurations. Therefore, our primary purpose is to show that our algorithms work and that our implementation is correct. 
This is demonstrated through qualitative comparisons and some quantitative results in terms of a
computational error analysis.

The outline of this work is as follows. In Section~\ref{sec2}, to start our discussion, we will briefly describe
the polynomials required for the N\'ed\'elec basis. 
Moreover, we give a short overview of current state-of-the-art methods of addressing the sign conflict on uniform grids. In Section~\ref{sec:extension_to_non-conforming}, we move on to non-conforming grids. 
In that section, we describe the modifications that are necessary to ensure global continuity even in the presence of hanging faces. We especially focus on the details required to implement a method to ensure global continuity.
Section~\ref{sec:constraint-matrix} is the key section of this work, describing the necessary modifications that have to be applied to the constraint matrix. 
We also provide a detailed explanation of how to overcome this sign conflict introduced from the constraint matrix, with some examples of pseudo-code.
Section~\ref{sec:numerical_tests} briefly introduces
the time-harmonic Maxwell's equations and substantiates our implementation with the help of two numerical examples.

%+-------------------------------------------------------------------+
%| Section: Preliminaries and Principal Problem of the Sign-Conflict |
%+-------------------------------------------------------------------+
\section{Preliminaries and Principal Problem of the Sign-Conflict}
\label{sec2}%
\subsection{$\mathbf{H}_\text{curl}$-conforming element space}
Let $\Omega\subset\mathbb{R}^d$, $d=2,3,$ be a bounded Lipschitz domain. The discretization of the Sobolev space \[
\mathbf{H}_\text{curl}(\Omega)=\{ \vec{u}\in [L^2(\Omega)]^d: \text{curl}\;\vec{u} \in [L^2(\Omega)]^{(2d-3)}\}, \quad d=2,3,
\]
requires tangential continuity along element interfaces.
The first and simplest conforming finite element spaces were developed by Nédélec  \cite{Art:Nedelec:80,Art:Nedelec:86}.
They preserve the tangential continuity. The systematic construction of higher-order FE spaces 
uses the De Rham cohomology. We refer the reader to \cite{Bk:Demkowicz.I:2008,Bk:Monk:2003} for more details.

For the polynomial basis, we choose Legendre~\cite{Bk:Szego:OrthogonalPolynomials:1939} and integrated
Legendre polynomials~\cite{Bk:Szabo:FiniteElementAnalysis:2021}, as they will
provide good sparsity properties in the involved element matrices \cite{Diss:Zaglmayr:2006}[Chapter 5.2.1]. For $n\geq 2$, we define the integrated Legendre polynomials by $L_n(x) \coloneqq \int_{-1}^{x} l_{n-1}(\xi) \mathsf{d} \xi$ for $x\in[-1,1]$,
where $l_p(x)=\frac{1}{2^p p!} \frac{\mathrm{d}^p}{\mathrm{d}x^p} (x^2-1)^p $ denotes the $p$-th Legendre polynomial.
Note that
\begin{equation}
  \begin{array}{rl}
    L_1(x)  &=  x + 1 , \\
    L_2(x)  &= \frac{1}{2}\left( x^2 - 1 \right), \\
    (n + 1) L_{n+1}(x)  &= (2n - 1)x L_n(x) - (n - 2) L_{n-1}(x), \quad\text{ for } n \geq 2, \quad x\in[-1,1].
  \end{array}
\end{equation}
This recursive formula allows an efficient point evaluation of the integrated Legendre polynomials. 
The concept of employing integrated Legendre polynomials as a polynomial basis for $\mathbf{H}_\text{curl}(\Omega)$ space was introduced in \cite{Art:Ainsworth:01} for quadrilateral elements.

\begin{figure}
  \begin{center}
    \input{Gnuplot/integrated_legendre.tex}
    \caption{%
      \label{fig:integrated_legendre}
      The integrated Legendre polynomials $L_2,~L_3,~L_4$ and $L_5$ are depicted. Integrated Legendre polynomials corresponding to even polynomial 
      degrees are symmetric, while those corresponding to odd polynomial degrees are antisymmetric.
    }
  \end{center}
\end{figure}

For three space dimensions, there are edge-, face- and cell-based basis functions. More precisely, the
cell-based basis functions on $\mathcal{C}^3=[0,1]^3$ up to the maximal polynomial degree $p_C$ are defined as
\begin{eqnarray*}
  \phi^{(curl,a)}_{i,j,k}(x_1,x_2,x_3)              & = &  \nabla_a(L_i(2x_1-1)\;L_j(2x_2-1)\;L_k(2x_3-1)),\\
  \phi^{(curl,IV)}_{i,j}(x_1,x_2,x_3)               & = & L_i(2x_\alpha) L_j(2x_\beta) \nabla x_\gamma,
\end{eqnarray*}
with $i,j,k=2,\ldots,p_C,~a \in\{I,II,III\},~(\alpha,\beta,\gamma)\in\{(1,2,3),(2,3,1),(3,2,1) \}$, the gradient $\nabla_I=\nabla$ and the antigradients $\nabla_{II}=\nabla_I-2\frac{\partial}{\partial x_2} (0,1,0)^\top$ and $\nabla_{III}=\nabla_I-2\frac{\partial}{\partial x_3} (0,0,1)^{\top}$.
In the same way, the other basis functions are defined. We refer to  
the work of Zaglmayr~\cite[Chapter 5.2]{Diss:Zaglmayr:2006} for a detailed definition.

\subsection{Reference Cell in Two Dimensions}
The enumeration of vertices and faces is based on the implementation in \dealii~\cite{Software:dealii:9.5}.
A more detailed description of the cell is given in the \dealii documentation 
\footnote{\url{https://www.dealii.org/current/doxygen/deal.II/structGeometryInfo.html}}.
We define the quadrilateral reference element as $\mathcal{C}^2 = [0, 1] \times [0, 1]$ with the default parametrization.
It is bounded by its faces.
As the vertex ordering is a crucial part of this work, we introduce the vertex enumeration on the 
reference cell in Figure \ref{fig:reference_cell}.
Moreover, we need the set of all faces, which is given by
$\mathcal{F} = \left\{ F_m \right\}_{0 \leq m < 4}$ 
with the local face-ordering
$F_m = \{ v_i, v_j \}$ where $(i,j) \in \{(0,2),~(1,3),~(0,1),~(2,3)\}$.
We denote the cell itself with local vertex-ordering $C = \{v_0, v_1, v_2, v_3\}$.
The polynomial degree vector is given by $\vec p = \left( \left\{ p_F \right\}_{F \in \mathcal{F}}, p_C \right)$.

\begin{figure}
  \begin{minipage}{0.45\textwidth}
    \begin{tikzpicture}[]
  % === Vertices ===
  % vertices
  \coordinate (V0) at (0.0, 0.0);
  \coordinate (V1) at (3.0, 0.0);
  \coordinate (V2) at (0.0, 3.0);
  \coordinate (V3) at (3.0, 3.0);

  % draw vertices:
  \foreach \n in {V0, V1, V2, V3}
    \node at (\n)[circle, fill, inner sep=1.5pt]{};

  % vertice labels
  \node[anchor=north west] at (V0) {$v_0$};
  \node[below] at (V1) {$v_1$};
  \node[anchor=south east] at (V2) {$v_2$};
  \node[above] at (V3) {$v_3$};

  % edge labels
  \node[anchor=east, gray]  at (-0.2,  1.5) {$F_0$};
  \node[anchor=west, gray]  at ( 3.2,  1.5) {$F_1$};
  \node[anchor=north, gray] at ( 1.5, -0.2) {$F_2$};
  \node[anchor=south, gray] at ( 1.5,  3.2) {$F_3$};

  % === Axis ===
  \coordinate (Vx) at (3.7, 0.0);
  \coordinate (Vy) at (0.0, 3.7);
  
  \draw[->, >=stealth] (-0.5, 0) -- (Vx);
  \draw[->, >=stealth] (0, -0.5) -- (Vy);
  
  \node[right] at (Vx) {x};
  \node[above] at (Vy) {y};

  % === Edges ===
  \draw[thick]      (V0) -- (V1);
  \draw[thick]      (V0) -- (V2);
  \draw[thick]      (V3) -- (V2);
  \draw[thick]      (V3) -- (V1);

  % === Direction of the edges ===
  \draw[->, >=stealth, thin, gray] (-0.2, 0.5) -- (-0.2, 2.5);
  \draw[->, >=stealth, thin, gray] (3.2, 0.5) -- (3.2, 2.5);
  \draw[->, >=stealth, thin, gray] (0.5, -0.2) -- (2.5, -0.2);
  \draw[->, >=stealth, thin, gray] (0.5, 3.2) -- (2.5, 3.2);

  % === Lightning Effect ===
  % face 0
  \filldraw[thin, fill=gray, fill opacity=0.25] (V0) -- (V1) -- (V3) -- (V2) -- cycle;
\end{tikzpicture}
  \end{minipage}
  \hfill
  \begin{minipage}{0.45\textwidth}
    \begin{tikzpicture}[]
  % === Vertices ===
  % vertices
  \coordinate (V0) at (0.0, 0.0);
  \coordinate (V1) at (3.0, 0.0);
  \coordinate (V2) at (1.0, 1.0);
  \coordinate (V3) at (4.0, 1.0);
  \coordinate (V4) at (0.0, 3.0);
  \coordinate (V5) at (3.0, 3.0);
  \coordinate (V6) at (1.0, 4.0);
  \coordinate (V7) at (4.0, 4.0);

  % edges
  \coordinate (E0) at (0.5, 0.5);
  \coordinate (E1) at (3.5, 0.5);
  \coordinate (E2) at (1.5, 0.0);
  \coordinate (E3) at (2.5, 1.0);
  \coordinate (E4) at (0.5, 3.5);
  \coordinate (E5) at (3.5, 3.5);
  \coordinate (E6) at (1.5, 3.0);
  \coordinate (E7) at (2.5, 4.0);
  \coordinate (E8) at (0.0, 1.5);
  \coordinate (E9) at (3.0, 1.5);
  \coordinate (E10) at (1.0, 2.5);
  \coordinate (E11) at (4.0, 2.5);

  % draw vertices:
  \foreach \n in {V0, V1, V3, V4, V5, V6, V7}
    \node at (\n)[circle, fill, inner sep=1.5pt]{};
  
  \node at (V2) [circle, fill, inner sep=1.5pt, gray]{};

  % === vertice lables ===
  \node[anchor=north west] at (V0) {$v_0$};
  \node[below] at (V1) {$v_1$};
  
  \node[below, gray, xshift=3pt] at (V2) {$v_2$};
  \node[below, xshift=3pt] at (V3) {$v_3$};
  
  \node[anchor=south east] at (V4) {$v_4$};
  \node[above, xshift=-5pt] at (V5) {$v_5$};
  
  \node[above] at (V6) {$v_6$};
  \node[above] at (V7) {$v_7$};

  % === Edge lables ===
  \node[right,        lightgray, yshift=-1pt] at (0.45, 0.3) {$E_0$};
  \node[right,        gray,      yshift=-1pt] at (3.45, 0.3) {$E_1$};
  \node[anchor=north, gray,      yshift=-3pt] at (E2) {$E_2$};
  \node[below,        lightgray, yshift=-3pt] at (E3) {$E_3$};

  \node[left,         gray, yshift=1pt] at (0.6, 3.8) {$E_4$};
  \node[left,         gray, yshift=1pt] at (3.55, 3.7) {$E_5$};
  \node[anchor=south, gray, yshift=3pt] at (E6) {$E_6$};
  \node[above,        gray, yshift=3pt] at (E7) {$E_7$};

  \node[left,  gray,      xshift=-3pt] at (E8) {$E_8$};
  \node[right, gray,      xshift= 3pt] at (E9) {$E_9$};
  \node[left,  lightgray, xshift=-3pt] at (E10) {$E_{10}$};
  \node[right, gray,      xshift= 3pt] at (E11) {$E_{11}$};

  % === Axis ===
  \coordinate (Vx) at (4.0, 0.0);
  \coordinate (Vy) at (1.5, 1.5);
  \coordinate (Vz) at (0.0, 4.0);
  
  \draw[->, >=stealth] (-0.5, 0) -> (Vx);
  \draw[->, thin, gray, >=stealth] (0, 0) -> (Vy);
  \draw[->, >=stealth] (0, -0.5) -> (Vz);
  \draw[thin] (0, 0) -> (-0.5, -0.5);
  
  \node[right] at (Vx) {x};
  \node[anchor=south west, gray] at (Vy) {y};
  \node[above] at (Vz) {z};
  
  % === Edges direction === 
  \draw[->, >=stealth, thin, gray] (-0.2, 0.5) -- (-0.2, 2.5);
  \draw[->, >=stealth, thin, gray] ( 3.2, 0.5) -- ( 3.2, 2.5);
  \draw[->, >=stealth, thin, gray] ( 0.5,-0.2) -- ( 2.5,-0.2);
  \draw[->, >=stealth, thin, gray] ( 0.5, 3.2) -- ( 2.5, 3.2);

  \draw[->, >=stealth, thin, lightgray] ( 0.8, 1.5) -- ( 0.8, 3.5);
  \draw[->, >=stealth, thin, gray] ( 4.2, 1.5) -- ( 4.2, 3.5);
  \draw[->, >=stealth, thin, lightgray] ( 1.5, 0.8) -- ( 3.5, 0.8);
  \draw[->, >=stealth, thin, gray] ( 1.5, 4.2) -- ( 3.5, 4.2);

  %\draw[->, >=stealth, thin, gray] ( 0.05, 0.25) -- ( 0.65, 0.85);
  \draw[->, >=stealth, thin, lightgray] ( 0.3, 0.1) -- ( 0.9, 0.7);
  \draw[->, >=stealth, thin, gray] ( 0.05, 3.25) -- ( 0.65, 3.85);
  \draw[->, >=stealth, thin, gray] ( 3.3, 0.1) -- ( 3.9, 0.7);
  %\draw[->, >=stealth, thin, gray] ( 3.3, 3.1) -- ( 3.9, 3.7);
  \draw[->, >=stealth, thin, gray] ( 3.05, 3.25) -- ( 3.65, 3.85);

  % === Edges ===
  % face 2
  \draw[thick]      (V0) -- (V1);
  \draw[thick]      (V0) -- (V4);
  \draw[thick]      (V5) -- (V4);
  \draw[thick]      (V5) -- (V1);
  
  % facethick,  3
  \draw[thin, gray] (V2) -- (V3);
  \draw[thin, gray] (V2) -- (V6);
  \draw[thick]      (V7) -- (V6);
  \draw[thick]      (V7) -- (V3);
  
  \draw[thin, gray] (V0) -- (V2);
  \draw[thick]      (V4) -- (V6);
  \draw[thick]      (V5) -- (V7);
  \draw[thick]      (V1) -- (V3);
  
  % === Lightning Effect ===
  % face 0
  \filldraw[thin, fill=gray, fill opacity=0.25] (V0) -- (V1) -- (V5) -- (V4) -- cycle;

  % face 3
  \filldraw[thin, fill=gray, fill opacity=0.13] (V1) -- (V3) -- (V7) -- (V5) -- cycle;

  % face 5
  \filldraw[thin, fill=gray, fill opacity=0.13] (V4) -- (V5) -- (V7) -- (V6) -- cycle;
  
  \node at (V0) [circle, fill, inner sep=1.5pt]{};
  \node at (V3) [circle, fill, inner sep=1.5pt]{};
  \node at (V6) [circle, fill, inner sep=1.5pt]{};

\end{tikzpicture}
  \end{minipage}
  \caption{%
    \label{fig:reference_cell}
    Left: Vertex and face ordering of the two-dimensional reference element.
    Right: Vertex, edge, and face ordering of the three-dimensional reference element.
  }
\end{figure}

\subsection{Reference Cell in Three Dimensions}
We define the reference element in three dimensions as
$\mathcal{C}^3 = [0, 1] \times [0, 1] \times [0,1]$ again with the default parametrization
and the vertex ordering shown in Figure \ref{fig:reference_cell}. 
The set of all edges is given by $\mathcal{E} = \left\{ E_m \right\}_{0 \leq m < 12}$
with local edge-ordering $E_m = \{ v_i, v_j \}$ as shown in Figure~\ref{fig:reference_cell}.
The local face order is given by
\begin{equation}
  \label{eq:faces}
  \begin{array}{r l}
    \mathcal{F} = \left\{ F_m \right\}_{0 \leq m < 6} =
    \{&\{v_0, v_2, v_4, v_6 \},
    \{v_1, v_3, v_5, v_7 \},
    \{v_0, v_1, v_4, v_5 \}, \\
    &\{v_2, v_3, v_6, v_7 \},
    \{v_0, v_1, v_2, v_3 \},
    \{v_4, v_5, v_6, v_7 \} ~ \}.
  \end{array}
\end{equation}
The polynomial degree vector is given by
$\vec p = \left( \left\{ p_E \right\}_{E \in \mathcal{E}}, \left\{ p_F \right\}_{F \in \mathcal{F}}, p_C \right)$.

\subsection{Principal Problem of the Sign Conflict}
\label{sec:principal_sign-conflict}
In this subsection, we briefly outline the fundamental idea behind the 
algorithm on how to overcome the sign conflict. Details will then be explained 
in the following sections.

To ensure the continuity between two neighboring elements, the resulting
polynomials on the faces in two dimensions and on the edges and the faces
in three dimensions must match. 
The integrated Legendre polynomials are either symmetric for even polynomial degrees or anti-symmetric for odd polynomial degrees; see Figure~\ref{fig:integrated_legendre}.
The FE map transforms the local basis functions of the reference element to one element of the mesh. All interior faces of a two-dimensional element share two neighboring elements.
Due to the FE map, a face with vertices $v_1$ and $v_2$ can either start from $v_2$ or from $v_1$. If some of the basis functions, as in our case, are not symmetric, the required
global continuity conditions of the global FE space would fail.

One solution to overcome the sign conflict
was proposed by Zaglmayr \cite{Diss:Zaglmayr:2006}
and implemented into \dealii by Kynch and Ledger
\cite{Art:Kynch:ResolvingTheSignConflict:2017}.
Their paper also provides some visualization of the sign conflict.
The basic idea of one possible algorithm that was proposed by Zaglmayr to solve the sign conflict on non-orientable grids is to use the global vertex indices to decide the orientation of edges and faces.
In any given mesh, each vertex is assigned to a unique global index by the finite element software.
When examining an edge or a face, these global vertex indices are taken into account.
For an edge or face of a two dimensional element, the two vertices are considered.
If the global index of the first vertex is smaller than that of the second, the orientation is done from the first vertex to the second.
Conversely, if the global index of the first vertex is larger, the orientation is done from the second vertex to the first.

For a face of a three-dimensional element, the direction of the outer lines is determined in a similar manner as for the edges.
However, one direction needs to be designated as the primary direction.
This is achieved by comparing the global vertex indices of the neighboring vertices of the first vertex of a face.
In Figure~\ref{fig:reference_cell}, this corresponds to $v_1$ and $v_2$.
If $v_1<v_2$, the $x$-direction is chosen as the primary direction.
If $v_1>v_2$, the $y$-direction is selected as the primary direction.
This approach ensures a consistent orientation across different elements, which is crucial for avoiding the sign conflict.

%+-------------------------------------------------------------------+
%|           Section: Global Continuity on Non-Conforming Grids      |
%+-------------------------------------------------------------------+
\section{Global Continuity on Non-Conforming Grids}
\label{sec:extension_to_non-conforming}

In this section, we explain an algorithm to ensure the global continuity of the N\'ed\'elec elements in the presence of hanging edges.
The basic idea was already provided in \cite{Art:Kynch:ResolvingTheSignConflict:2017}.
Hence, we focus mainly on the essential details of the implementation. 
Moreover, we introduce Algorithm \ref{alg:find_hidden_hanging_edges} to cover all special cases as well.

\subsection{Identification of Hanging Faces}
We split the task of ensuring global continuity into two subproblems.
First, we identify all hanging faces and edges, and later, in Subsection \ref{subsec:modify_hanging_edges_and_faces}, we discuss how to modify those hanging edges and faces in order to ensure global continuity.
A face $F$ is called a hanging face if and only if the neighboring face $N_F$ is coarser than $F$.
To identify all hanging faces, we loop over all cells $K$ in the grid $\mathcal{K}$ and mark all faces that have a coarser neighbor as hanging faces.

\subsection{Identification of Hanging Edges in Three Dimensions}
In the three-dimensional case, we also have to consider hanging edges. 
Here, the definition is similar: an edge $E$ is called a hanging edge if and only if the neighboring edge $N_E$ is coarser than $E$.

\begin{figure}
  \begin{center}
    \includegraphics[width=0.5\textwidth]{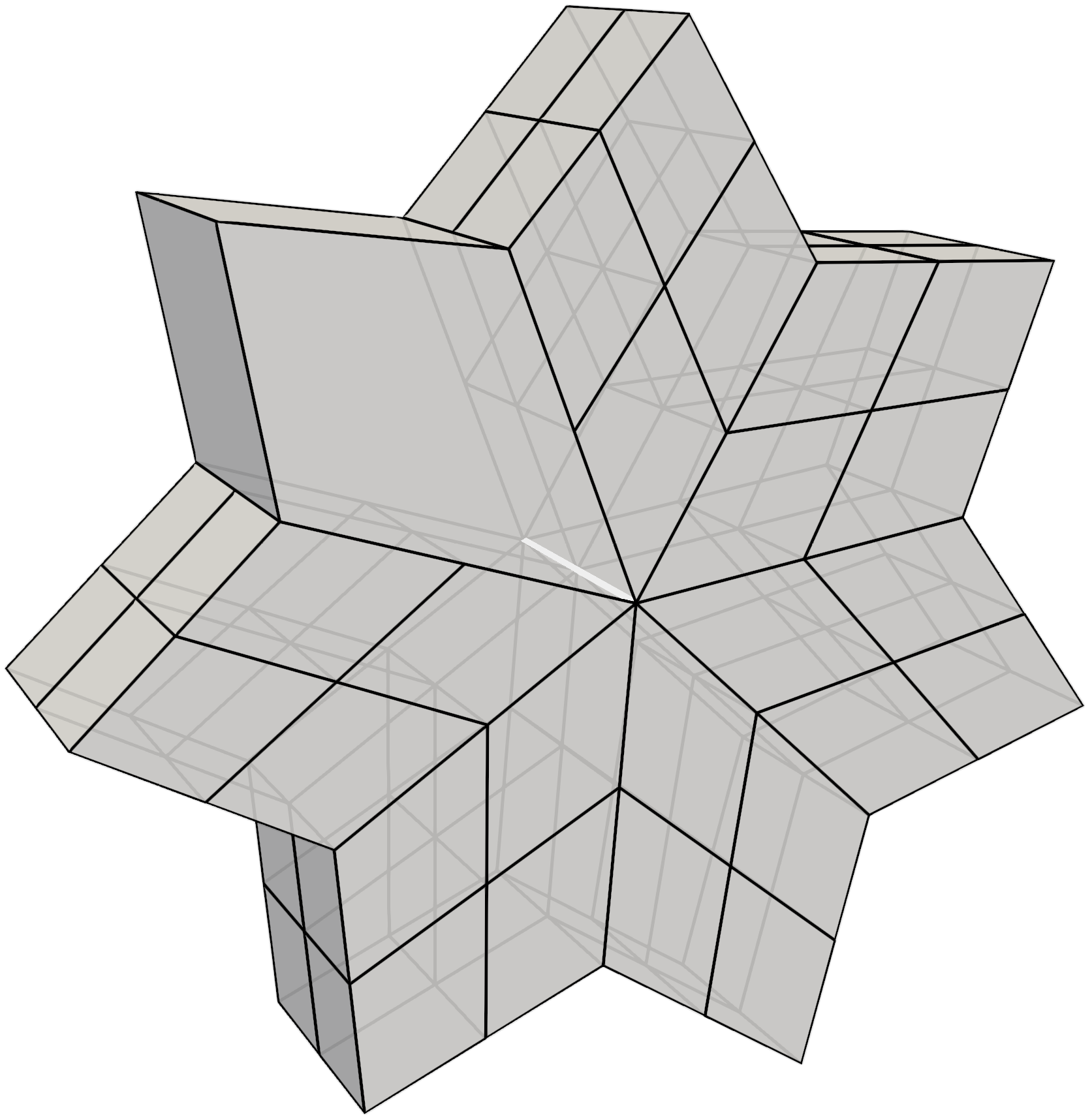}
  \end{center}
  \caption{
    Most cells have no hanging faces but a hanging edge.
    \label{fig:hanging_edge}
  }
\end{figure}

In three dimensions, certain configurations may result in an element having an edge that neighbors a coarser element,
even though the neighbors of all faces of that element are of the same refinement level.
This can lead to the presence of hanging edges that do not belong to a hanging face.
An example of such a configuration, where seven cells share a common edge, is shown in Figure~\ref{fig:hanging_edge}.
The algorithm to find these hanging edges is presented in Algorithm~\ref{alg:find_hidden_hanging_edges}.

\begin{algorithm}
  \caption{\label{alg:find_hidden_hanging_edges}Find remaining hanging edges}
  \LoopOverAll ( cells $K$ in grid $\mathcal{K}$ \textbf{do} ) {
    \LoopOverAll ( edges $E\in\mathcal{E}$ from the current cell $K$ \textbf{do} ) {
      Skip all edges that do belong to a hanging face; \\
      \LoopOverAll ( neighbour cells $N\in \mathcal{N}_E$ that are adjacent to the current edge $E$ \textbf{do} ) {
        \If { The neighbour cell $N$ is coarser than the current cell $K$ } {
          Mark the edge $E$ as the hanging edge.
        }
      }
    }
  }
\end{algorithm}

\subsection{Adapting Cell Orientation in the Presence of Hanging Faces and Edges}
\label{subsec:modify_hanging_edges_and_faces}
After identifying all hanging faces and edges, it is crucial to adapt their orientation to ensure the continuity of the mesh.
This process is outlined in Algorithm~\ref{alg:hanging_edge_direction}.
Figure~\ref{fig:sign_conflict_2d} illustrates the difference in grid orientation with and without this special treatment for hanging edges.

\begin{figure}[h]
  \begin{minipage}{0.45\textwidth}
    \begin{center}
      \begin{tikzpicture}[]
  \draw[thick] (0.00, 0.00) -- (4.00, 0.00);
  \draw[thick] (0.00, 2.00) -- (4.00, 2.00);
  \draw[thick] (0.00, 4.00) -- (4.00, 4.00);
  \draw[thick] (0.00, 0.00) -- (0.00, 4.00);
  \draw[thick] (2.00, 0.00) -- (2.00, 4.00);
  \draw[thick] (4.00, 0.00) -- (4.00, 4.00);

  \draw[thick] (0.00, 1.00) -- (2.00, 1.00);
  \draw[thick] (1.00, 0.00) -- (1.00, 2.00);

  % Directions:
  % bottom
  \draw[->, >=stealth] (0.25, 0.10) -- (0.75, 0.10);
  \draw[<-, >=stealth] (1.25, 0.10) -- (1.75, 0.10);
  \draw[<-, >=stealth] (2.25, 0.10) -- (3.75, 0.10);

  % Refined
  \draw[->, >=stealth] (0.25, 0.90) -- (0.75, 0.90);
  \draw[<-, >=stealth] (1.25, 0.90) -- (1.75, 0.90);
  \draw[->, >=stealth] (0.25, 1.10) -- (0.75, 1.10);
  \draw[<-, >=stealth] (1.25, 1.10) -- (1.75, 1.10);

  % center lower
  \draw[->, >=stealth] (0.25, 1.90) -- (0.75, 1.90);
  \draw[<-, >=stealth] (1.25, 1.90) -- (1.75, 1.90);
  \draw[<-, >=stealth] (2.25, 1.90) -- (3.75, 1.90);

  % center higher
  \draw[->, >=stealth] (0.25, 2.10) -- (1.75, 2.10);
  \draw[<-, >=stealth] (2.25, 2.10) -- (3.75, 2.10);

  % top
  \draw[->, >=stealth] (0.25, 3.90) -- (1.75, 3.90);
  \draw[<-, >=stealth] (2.25, 3.90) -- (3.75, 3.90);

  % Left
  \draw[->, >=stealth] (0.10, 0.25) -- (0.10, 0.75);
  \draw[<-, >=stealth] (0.10, 1.25) -- (0.10, 1.75);
  \draw[<-, >=stealth] (0.10, 2.25) -- (0.10, 3.75);

  % Refined
  \draw[->, >=stealth] (0.90, 0.25) -- (0.90, 0.75);
  \draw[<-, >=stealth] (0.90, 1.25) -- (0.90, 1.75);
  \draw[->, >=stealth] (1.10, 0.25) -- (1.10, 0.75);
  \draw[<-, >=stealth] (1.10, 1.25) -- (1.10, 1.75);

  % center - left
  \draw[->, >=stealth] (1.90, 0.25) -- (1.90, 0.75);
  \draw[<-, >=stealth] (1.90, 1.25) -- (1.90, 1.75);
  \draw[<-, >=stealth] (1.90, 2.25) -- (1.90, 3.75);

  % center - right
  \draw[->, >=stealth] (2.10, 0.25) -- (2.10, 1.75);
  \draw[<-, >=stealth] (2.10, 2.25) -- (2.10, 3.75);

  % right
  \draw[->, >=stealth] (3.90, 0.25) -- (3.90, 1.75);
  \draw[<-, >=stealth] (3.90, 2.25) -- (3.90, 3.75);
\end{tikzpicture}
    \end{center}
  \end{minipage}
  \hfill
  \begin{minipage}{0.45\textwidth}
    \begin{center}
      \begin{tikzpicture}[]
  \draw[thick] (0.00, 0.00) -- (4.00, 0.00);
  \draw[thick] (0.00, 2.00) -- (4.00, 2.00);
  \draw[thick] (0.00, 4.00) -- (4.00, 4.00);
  \draw[thick] (0.00, 0.00) -- (0.00, 4.00);
  \draw[thick] (2.00, 0.00) -- (2.00, 4.00);
  \draw[thick] (4.00, 0.00) -- (4.00, 4.00);

  \draw[thick] (0.00, 1.00) -- (2.00, 1.00);
  \draw[thick] (1.00, 0.00) -- (1.00, 2.00);

  % Directions:
  % bottom
  \draw[->, >=stealth] (0.25, 0.10) -- (0.75, 0.10);
  \draw[<-, >=stealth] (1.25, 0.10) -- (1.75, 0.10);
  \draw[<-, >=stealth] (2.25, 0.10) -- (3.75, 0.10);

  % Refined
  \draw[->, >=stealth] (0.25, 0.90) -- (0.75, 0.90);
  \draw[<-, >=stealth] (1.25, 0.90) -- (1.75, 0.90);
  \draw[->, >=stealth] (0.25, 1.10) -- (0.75, 1.10);
  \draw[<-, >=stealth] (1.25, 1.10) -- (1.75, 1.10);

  % center lower
  \draw[->, >=stealth] (0.25, 1.90) -- (0.75, 1.90);
  \draw[->, >=stealth] (1.25, 1.90) -- (1.75, 1.90);
  \draw[<-, >=stealth] (2.25, 1.90) -- (3.75, 1.90);

  % center higher
  \draw[->, >=stealth] (0.25, 2.10) -- (1.75, 2.10);
  \draw[<-, >=stealth] (2.25, 2.10) -- (3.75, 2.10);

  % top
  \draw[->, >=stealth] (0.25, 3.90) -- (1.75, 3.90);
  \draw[<-, >=stealth] (2.25, 3.90) -- (3.75, 3.90);

  % Left
  \draw[->, >=stealth] (0.10, 0.25) -- (0.10, 0.75);
  \draw[<-, >=stealth] (0.10, 1.25) -- (0.10, 1.75);
  \draw[<-, >=stealth] (0.10, 2.25) -- (0.10, 3.75);

  % Refined
  \draw[->, >=stealth] (0.90, 0.25) -- (0.90, 0.75);
  \draw[<-, >=stealth] (0.90, 1.25) -- (0.90, 1.75);
  \draw[->, >=stealth] (1.10, 0.25) -- (1.10, 0.75);
  \draw[<-, >=stealth] (1.10, 1.25) -- (1.10, 1.75);

  % center - left
  \draw[->, >=stealth] (1.90, 0.25) -- (1.90, 0.75);
  \draw[->, >=stealth] (1.90, 1.25) -- (1.90, 1.75);
  \draw[<-, >=stealth] (1.90, 2.25) -- (1.90, 3.75);

  % center - right
  \draw[->, >=stealth] (2.10, 0.25) -- (2.10, 1.75);
  \draw[<-, >=stealth] (2.10, 2.25) -- (2.10, 3.75);

  % right
  \draw[->, >=stealth] (3.90, 0.25) -- (3.90, 1.75);
  \draw[<-, >=stealth] (3.90, 2.25) -- (3.90, 3.75);
\end{tikzpicture}
    \end{center}
  \end{minipage}
  \caption{
    \label{fig:sign_conflict_2d}
    Comparison of grid orientations.
    The left-hand side shows the grid without special treatment for hanging edges,
    while the right-hand side shows the grid with special treatment for hanging edges.
    }
\end{figure}

\begin{algorithm}
  \caption{\label{alg:hanging_edge_direction}Adapt the cell orientation in the presence of hanging faces and edges.}
  \LoopOverAll ( cells $K$ in grid $\mathcal{K}$ \textbf{do} ) {
    \LoopOverAll ( faces $F\in\mathcal{F}$ from cell $K$ \textbf{do} ) {
      \eIf{ face $F$ is marked as hanging face } {
        Compute the face orientation based on the global vertex indices of the parent cell of cell $K$;
      } {
        \tcp{ face $F$ is not marked as hanging face }
        Compute the face orientation based on the global vertex indices of cell $K$;
      }
    }
    \If{dim == 3}{
      \LoopOverAll ( edges of $E\in\mathcal{E}$ from cell $K$ ) {
        \eIf{ edge $E$ is marked as hanging edge } {
          Compute the edge orientation based on the global vertex indices of the parent cell of cell $K$;
        } {
          \tcp{ edge $E$ is not marked as hanging edge }
          Compute the edge orientation based on the global vertex indices of cell $K$;
        }
      }
    }
  }
\end{algorithm}

%+-------------------------------------------------------------------+
%|          Section: Modifications of the constraint matrix          |
%+-------------------------------------------------------------------+
\section{Modifications of the constraint matrix}
\label{sec:constraint-matrix}

\subsection{Solving the Mismatch Between the Number of Degrees of Freedom of Refined and Coarse Elements}
When a structured mesh is locally refined, hanging faces are introduced, and in the three-dimensional case, hanging edges are introduced as well.
This leads to a mismatch between the number of degrees of freedom (DoFs) of the refined and coarse elements.
The most prominent approach to deal with these additional DoFs is to impose constraints on the additional DoFs of the refined element by expressing them as a linear combination of the coarse element's DoFs. This can be written as
\begin{equation}\label{eq:constr}
  \varphi_r = [\alpha_{i,j}]_{i,j}^{n,m} \cdot \varphi_c,
\end{equation}
with $r$ denoting `refined', $c$ denoting `coarse', $i,j$ denoting the Dof indices,
$n$ and $m$ are the number of local DoFs involved in the constraints.
In more detail, $\varphi_r$ is the vector of the basis function on the refined element, $\varphi_c$ is the vector of the basis functions on
the coarse element, and $\alpha_{i,j}$ is the constraint matrix containing the weights between the corresponding basis functions.
The computation of the weights is not within the scope of this work for which we refer the reader to \cite{Art:Haubold:RecuurencesHighOrder:22, Art:Kus:14, Art:Stolfo:16}.
In the following, let us assume that the generation of the entries $\alpha_{i,j}$ is performed by a subroutine called \textbf{get\_local\_constraint\_matrix},
corresponding to the \dealii function FETools::compute\_face\_embedding\_matrices()\footnote{\url{https://www.dealii.org/current/doxygen/deal.II/namespaceFETools.html\#ac0fe5c7f55db091a4477af7c3989b83c}}.
Moreover, \textbf{local\_constraint\_to\_global} distributes the local constraint matrix into the global constraint matrix, 
which is more complicated in practice. Mostly, this corresponds to the \dealii function AffineConstraints::add\_entry()\footnote{\url{https://www.dealii.org/current/doxygen/deal.II/classAffineConstraints.html\#a2b7756e9cb8e53553211add5426f8e50}}.

When applying algorithms to ensure global continuity, we inevitably modify the orientation of the cells.
However, as discussed in Section~\ref{sec:extension_to_non-conforming}, 
the constraint matrix is computed
for the canonical coarse-fine mapping, which assumes a specific cell orientation.
Therefore, we must adjust the constraint matrix to respect the cells' 
orientation to prevent sign conflicts.
This adjustment involves multiplying the correct entries of the constraint matrix by $-1$.
In most cases, it is sufficient to compare the orientation of the parent cell to the child cell
rather than with the canonical orientation used in the canonical coarse-fine mapping.
If both the parent and the child cell differ from the canonical orientation,
the sign changes cancel each other out.
Therefore, it only matters whether the orientation of the parent and the child matches.

\begin{figure}
  \begin{minipage}{0.45\textwidth}
    \includegraphics[width=1.0\textwidth]{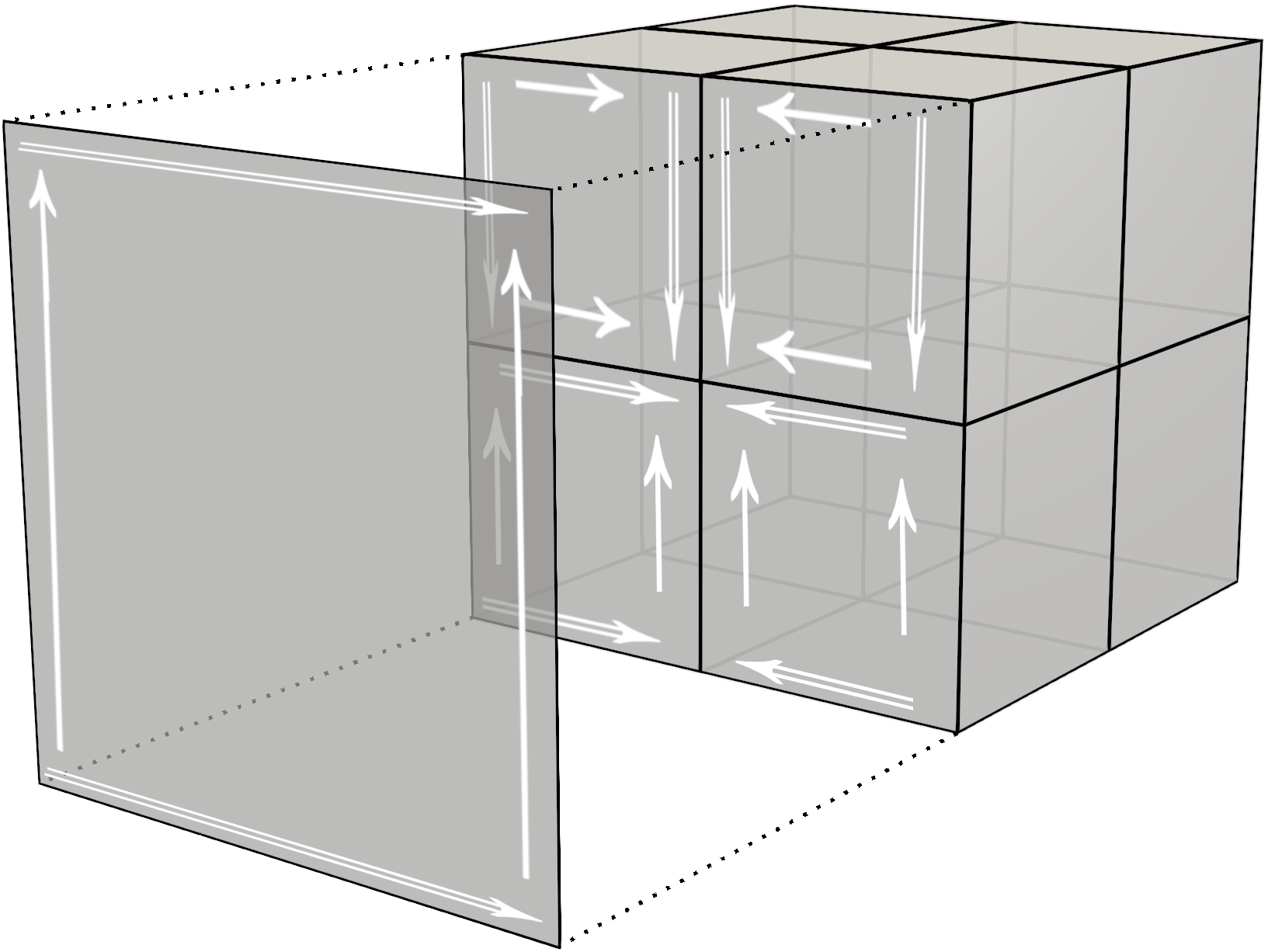}
  \end{minipage}
  \hfill
  \begin{minipage}{0.45\textwidth}
    \includegraphics[width=1.0\textwidth]{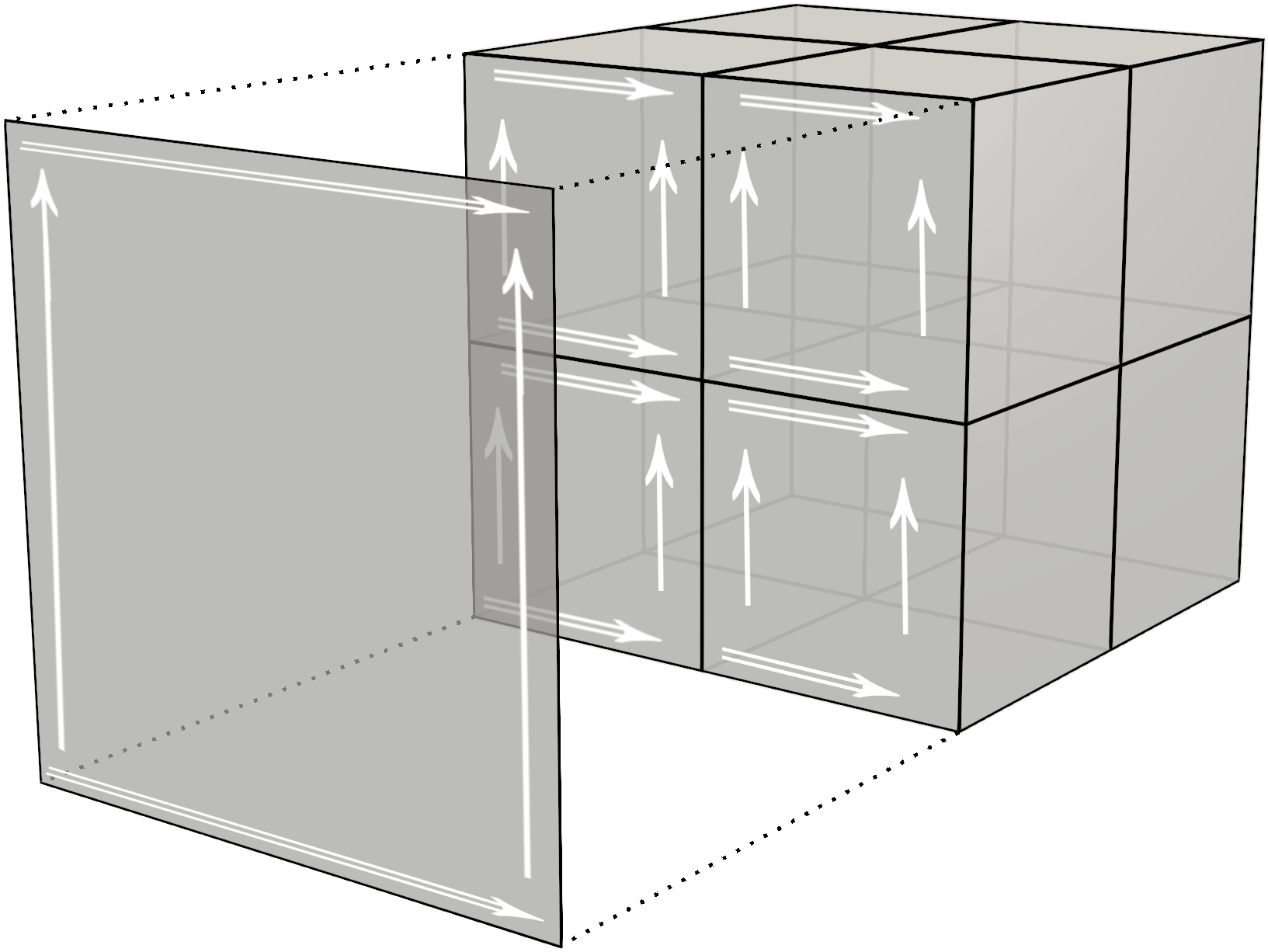}
  \end{minipage}
  \caption{%
    \label{fig:coarse_fine_configuartion}
    Left: Natural coarse-fine mapping, resulting from ignoring the hanging edges.
    Right: Canonical coarse-fine mapping.
  }
\end{figure}

Section~\ref{sec:principal_sign-conflict} summarized how to modify the grid to ensure global continuity.
In the case of the N\'ed\'elec elements, ensuring continuity on non-conforming meshes requires that the tangential components of the
basis function on the hanging edges and faces match those of the corresponding basis functions on the neighboring unrefined element.
The constraint matrix can be developed by considering a reference setting where we match the tangential constraints.
This reference setting is called the canonical-coarse mapping; see Figure~\ref{fig:coarse_fine_configuartion}.

These constraints can be applied to more general shapes with the help of an affine coordinate transformation \cite{Bk:Demkowicz.II:2008}.
The implementation presented in this work\footnote{\url{https://zenodo.org/records/10913219}} 
was created using \dealii as a programming platform that provides the 
functionality to
compute the weights numerically.
Therefore, we focus on modifying the given weights to match the grid's orientation.
The constraints for the hanging edges and faces
depend on the orientation of the refined element and its unrefined neighbors.
Consequently, the constraints have to be computed during the runtime of the numerical simulations.

\subsection{Constraints for Hanging Faces in Two Dimensions}
To implement the actual cell orientation, 
we begin by considering the faces of two-dimensional elements.
We compare the vertex order of the refined element with that of the coarse neighbor, similar to Algorithm \ref{alg:hanging_edge_direction}.
If the vertex order between the refined and coarse neighbors does not match, we must adapt the constraint matrix accordingly.
First of all, the local and global indices of the involved basis functions are required. 
The local indices $i$ and $j$ in \eqref{eq:constr} correspond either to symmetric or antisymmetric basis functions. 
The entry $\alpha_{ij}$ has to be multiplied with $-1$ if the pair $(i,j)$ corresponds to a pair of symmetric and antisymmetric basis functions. 
If both are antisymmetric, then this operation has to be done twice and cancels out.

\begin{figure}[h]
  \begin{minipage}{0.30\textwidth}
    \begin{displaymath}
      \left(\begin{array}{c}
              \varphi^{F_0^R}_1 \\[2mm]
              \varphi^{F_0^R}_{2}             \\[2mm]
              \varphi^{F_1^R}_1 \\[2mm]
              \varphi^{F_1^R}_{2}
      \end{array}\right)
      =
      \left(\begin{array}{cc}
              \sfrac{1}{2} &  \sfrac{1}{2} \\[2mm]
              0            &  \sfrac{1}{4} \\[2mm]
              \sfrac{1}{2} &  \sfrac{1}{2} \\[2mm]
              0            &  \sfrac{1}{4}
      \end{array}\right)
      \cdot
      \left(\begin{array}{c}
              \varphi^{F^C}_1 \\[2mm]
              \varphi^{F^C}_{2}
      \end{array}\right)
    \end{displaymath}
  \end{minipage}
  \hfill
  \begin{minipage}{0.15\textwidth}
    \begin{center}
      \begin{tikzpicture}[rotate=90]
  % Children
  \draw[thick, ->, >=stealth] ( 0.20,  0.075) -- (1.30,  0.075);
  \draw[thick, ->, >=stealth] ( 1.70,  0.075) -- (2.80,  0.075);

  \node at (1.5,0.075) [circle, fill, inner sep=1.2pt]{};
  \node[anchor=east, xshift=-2pt] at (1.5, 0.0) {$v_2$};

  \node[anchor=east] at (0.75, 0.0) {$F_0^R$};
  \node[anchor=east] at (2.25, 0.0) {$F_1^R$};

  % Parent
  \draw[thick, ->, >=stealth] ( 0.2, -0.075) -- (2.8, -0.075);

  \node at (0,0) [circle, fill, inner sep=2.0pt] {};
  \node at (3,0) [circle, fill, inner sep=2.0pt] {};

  \node[anchor=west] at (0.0, 0.0) {$v_0$};
  \node[anchor=west] at (3.0, 0.0) {$v_1$};

  \node[anchor=west] at (1.5, 0.0) {$F^C$};
\end{tikzpicture}
    \end{center}
  \end{minipage}
  \hfill
  \begin{minipage}{0.30\textwidth}
    \begin{displaymath}
      \left(\begin{array}{c}
              \varphi_1^{\widetilde F_0^R} \\[2mm]
              \varphi^{\widetilde F_0^R}_{2}             \\[2mm]
              \varphi_1^{\widetilde F_1^R} \\[2mm]
              \varphi^{\widetilde F_1^R}_{2}
      \end{array}\right)
      =
      \left(\begin{array}{cc}
              \sfrac{1}{2} &  \sfrac{1}{2} \\[2mm]
              0            &  \sfrac{1}{4} \\[2mm]
              \sfrac{1}{2} & -\sfrac{1}{2} \\[2mm]
              0            &  \sfrac{1}{4}
      \end{array}\right)
      \cdot
      \left(\begin{array}{c}
              \varphi^{\widetilde F^C}_1 \\[2mm]
              \varphi^{\widetilde F^C}_{2}
      \end{array}\right)
    \end{displaymath}
  \end{minipage}
  \hfill
  \begin{minipage}{0.15\textwidth}
    \begin{center}
      \begin{tikzpicture}[rotate=90]
  % Children
  \draw[thick, ->, >=stealth] ( 0.20,  0.075) -- (1.30,  0.075);
  \draw[thick, <-, >=stealth] ( 1.70,  0.075) -- (2.80,  0.075);

  \node at (1.5, 0.075) [circle, fill, inner sep=1.2pt]{};
  \node[anchor=east, xshift=-2pt] at (1.5, 0.0) {$\tilde v_2$};

  \node[anchor=east] at (0.75, 0.0) {$\widetilde F_0^R$};
  \node[anchor=east] at (2.25, 0.0) {$\widetilde F_1^R$};

  % Parent
  \draw[thick, ->, >=stealth] ( 0.2, -0.075) -- (2.8, -0.075);

  \node at (0,0) [circle, fill, inner sep=2.0pt] {};
  \node at (3,0) [circle, fill, inner sep=2.0pt] {};

  \node[anchor=west] at (0.0, 0.0) {$\tilde v_0$};
  \node[anchor=west] at (3.0, 0.0) {$\tilde v_1$};

  \node[anchor=west] at (1.5, 0.0) {$\widetilde F^C$};
\end{tikzpicture}
    \end{center}
  \end{minipage}
  \caption{
    \label{fig:constraint_matrix_2d_p2}
    The resulting constraint matrix, when Algorithm~\ref{alg:cm_2dim} is applied to
    a hanging face of a two-dimensional element with polynomial degree $p_F=2$.
    Note that $\phi_i^F$ is the basis function of face $F$ with degree $i$, $i=1,2$.
        Left: Canonical orientation. Right: One face differs from the canonical orientation.
  }
\end{figure}

Based on this information, we can formulate Algorithm~\ref{alg:cm_2dim} to resolve the sign conflict on hanging edges.
Note that here, we focus on one hanging face.
In an actual real-world implementation, one would need to loop over all faces and check if each face is marked as a hanging face.
Since hanging node constraints are only necessary for hanging faces, we assume that the outer loop for identifying hanging faces
is already implemented, and we concentrate on the inner part.
\begin{algorithm}
  \caption{
    \label{alg:cm_2dim}
    Given a face $F$ that is marked as a hanging face, adapt the constraint matrix based on the orientation of the face
    $F$ and the orientation of the children $F^R$ of face $F$.
  }
  \LoopOverAll ( children $F^R$ of face $F$ \textbf{do}) {
    \If{ The orientation of $F^R$ and $F$ does not match } {
      \tcp{Get the part of the constraint matrix that corresponds to the child $F^R$}
      local\_constraint\_matrix $\gets$ \textbf{get\_local\_constraint\_matrix}($F^R$)\;
      \tcp{Modify all constraint matrix entries that belong to this face and to anti-symmetric shape functions}
      \For{$i,j$ in local\_constraint\_matrix } {
        \If{ \textbf{is\_odd}(i + j) } {
          local\_constraint\_matrix(i, j) $\gets$ - local\_constraint\_matrix(i, j);
        }
      }
      \tcp{Write the modified local sub-constraint matrix into the global constraint matrix}
      \textbf{local\_constraint\_to\_global}(local\_constraint\_matrix);
    } 
  }
\end{algorithm}

\subsection{Constraints for Hanging Faces in Three Dimensions}
In our previous discussion, we focused solely on the orientation of hanging faces in two dimensions,
which corresponds to the edges in three dimensions. These hanging faces consist of eight external edges,
four internal edges, and four faces.
The face of the coarse element, on the other hand, consists of four external 
edges and one face. Consequently,
the size of the constraint matrix increases accordingly.

As the constraint matrix grows significantly in size for hanging faces, especially in the first
non-trivial case where the polynomial degree is $p_F=2$, we will only visualize the structure of
the constraint matrix in Figure~\ref{fig:constraint_matrix_structure}.

\begin{figure}
 \begin{minipage}{0.45\textwidth}
   \begin{center}
     \begin{tikzpicture}[]
  \draw[thick]               ( 0.1, 0.0) -- (2.9, 0);
  \draw[thick]               ( 0.1, 3.0) -- (2.9, 3);
  \draw[thick]               ( 0.0, 0.1) -- (0.0, 2.9);
  \draw[thick]               ( 3.0, 0.1) -- (3.0, 2.9);

  \node at (0,0) [circle, fill, inner sep=1.5pt]{};
  \node at (0,3) [circle, fill, inner sep=1.5pt]{};
  \node at (3,0) [circle, fill, inner sep=1.5pt]{};
  \node at (3,3) [circle, fill, inner sep=1.5pt]{};

  \node[anchor=north] at (1.5, 0.0) {$E_0^C$};
  \node[anchor=south] at (1.5, 3.0) {$E_1^C$};
  \node[anchor=east ] at (0.0, 1.5) {$E_2^C$};
  \node[anchor=west ] at (3.0, 1.5) {$E_3^C$};

  \node[            ] at (1.5, 1.5) {$F^C$};
\end{tikzpicture}
   \end{center}
 \end{minipage}
 \begin{minipage}{0.45\textwidth}
   \begin{center}
     \begin{tikzpicture}[]
  \draw[thick]               ( 0.1, 0.0) -- (1.4, 0);
  \draw[thick]               ( 1.6, 0.0) -- (2.9, 0);
  \draw[thick]               ( 0.1, 3.0) -- (1.4, 3);
  \draw[thick]               ( 1.6, 3.0) -- (2.9, 3);
  \draw[thick]               ( 0.0, 0.1) -- (0.0, 1.4);
  \draw[thick]               ( 0.0, 1.6) -- (0.0, 2.9);
  \draw[thick]               ( 3.0, 0.1) -- (3.0, 1.4);
  \draw[thick]               ( 3.0, 1.6) -- (3.0, 2.9);

  \draw[thick]               ( 0.1, 1.5) -- (1.4, 1.5);
  \draw[thick]               ( 1.6, 1.5) -- (2.9, 1.5);
  \draw[thick]               ( 1.5, 0.1) -- (1.5, 1.4);
  \draw[thick]               ( 1.5, 1.6) -- (1.5, 2.9);

  \node at (0,0) [circle, fill, inner sep=1.5pt]{};
  \node at (0,3) [circle, fill, inner sep=1.5pt]{};
  \node at (3,0) [circle, fill, inner sep=1.5pt]{};
  \node at (3,3) [circle, fill, inner sep=1.5pt]{};

  \node at (0.0, 1.5) [circle, fill, inner sep=1.5pt]{};
  \node at (3.0, 1.5) [circle, fill, inner sep=1.5pt]{};
  \node at (1.5, 0.0) [circle, fill, inner sep=1.5pt]{};
  \node at (1.5, 3.0) [circle, fill, inner sep=1.5pt]{};
  \node at (1.5, 1.5) [circle, fill, inner sep=1.5pt]{};

  \node[anchor=north] at (0.75, 0.00) {\small{$E^R_4$}};
  \node[anchor=north] at (2.25, 0.00) {\small{$E^R_5$}};
  \node[anchor=south] at (0.75, 3.00) {\small{$E^R_6$}};
  \node[anchor=south] at (2.25, 3.00) {\small{$E^R_7$}};
  \node[anchor=east ] at (0.00, 0.75) {\small{$E^R_8$}};
  \node[anchor=east ] at (0.00, 2.25) {\small{$E^R_9$}};
  \node[anchor=west ] at (3.00, 0.75) {\small{$E^R_{10}$}};
  \node[anchor=west ] at (3.00, 2.25) {\small{$E^R_{11}$}};

  \node[anchor=north, yshift=3pt] at (0.75, 1.50) {\small{$E^R_0$}};
  \node[anchor=north, yshift=3pt] at (2.25, 1.50) {\small{$E^R_1$}};
  \node[anchor=east , xshift=3pt] at (1.50, 0.75) {\small{$E^R_2$}};
  \node[anchor=east , xshift=3pt] at (1.50, 2.25) {\small{$E^R_3$}};

  \node[] at (0.65, 0.65) {\small{$F^R_0$}};
  \node[] at (2.25, 0.65) {\small{$F^R_1$}};
  \node[] at (0.65, 2.25) {\small{$F^R_2$}};
  \node[] at (2.25, 2.25) {\small{$F^R_3$}};
\end{tikzpicture}
   \end{center}
 \end{minipage}
 \caption{
   \label{fig:enumeration_of_edges}
   On the left-hand side is the enumeration of edges of the coarse parent face. On the right-hand side is the enumeration of edges and faces of the refined child faces.
 }
\end{figure}

\renewcommand{\arraystretch}{1.3}
\begin{figure}
  \begin{displaymath}
    \left(\begin{array} {c}
            l(E^R_0)     \\
            l(E^R_1)     \\
            l(E^R_2)     \\
            l(E^R_3)     \\
            l(E^R_4)     \\
            \vdots       \\
            l(E^R_7)     \\
            l(E^R_8)     \\
            \vdots       \\
            l(E^R_{11})  \\
            l(F^R_0)     \\
            \vdots       \\
            l(F^R_3)
    \end{array}
    \right)
    =
    \left(\begin{array} {c c c c c}
            C_{(0,0)} & C_{(0,1)} & C_{(0,2)}  & C_{(0,3)}  & C_{(0,4)} \\
            C_{(1,0)} & C_{(1,1)} & C_{(1,2)}  & C_{(1,3)}  & C_{(1,4)} \\
            C_{(2,0)} & C_{(2,1)} & C_{(2,2)}  & C_{(2,3)}  & C_{(2,4)} \\
            C_{(3,0)} & C_{(3,1)} & C_{(3,2)}  & C_{(3,3)}  & C_{(3,4)} \\
            C_{(4,0)} & C_{(4,1)} & 0          & 0          & 0         \\
            \vdots    & \vdots    & \vdots     & \vdots     & \vdots    \\
            C_{(7,0)} & C_{(7,1)} & 0          & 0          & 0         \\
            0         & 0         & C_{(8,2)}  & C_{(8,3)}  & 0         \\
            \vdots    & \vdots    & \vdots     & \vdots     & \vdots    \\
            0         & 0         & C_{(11,2)} & C_{(11,3)} & 0         \\
            0         & 0         & 0          & 0          & C_{(12,4)}\\
            \vdots    & \vdots    & \vdots     & \vdots     & \vdots    \\
            0         & 0         & 0          & 0          & C_{(15,4)}
    \end{array}
    \right)
    \cdot
    \left(\begin{array} {c}
            r(E^C_0) \\
            r(E^C_1) \\
            r(E^C_2) \\
            r(E^C_3) \\
            r(F^C)   \\
    \end{array}\right)
  \end{displaymath}
  \caption{
    \label{fig:constraint_matrix_structure}
    The structure of the constraint matrix. As a simplification, we group the basis functions as follows:
    $l(E^R_i)$, where $i \in \{0,11\}$, denotes the vector of all basis functions corresponding to the edge $E^R_i$ on the refined element.
    Next, $l(F^R_i)$, where $i \in \{0, 3\}$, denotes the vector of all basis functions on the face $F^R_i$.
    Similarly, $l(E^C_i)$, where $i \in \{0,3\}$, denotes the vector of all basis functions corresponding to the edge $E^C_i$ on the coarse element.
    Then, $l(F^R_i)$, where $i \in \{0, 3\}$, denotes the vector of all basis functions on the face $F^R_i$.
    Finally, $C_{(i,j)}$ represents the corresponding sub-constraint matrix between $l(E^C_i)$ and $r(E^R_j)$.
    The notation follows Figure~\ref{fig:enumeration_of_edges}.
}
\end{figure}
\renewcommand{\arraystretch}{1.0}

\subsection{Resolving the Sign Conflict on Hanging Faces in Three Dimensions}
Due to the complexity of the constraint matrix structure, we consider the different sub-constraint matrices,
denoted as $C_{(i,j)}$ in Figure \ref{fig:constraint_matrix_structure}, independently.
For each hanging edge and face, we determine which coarse edge and face directions must be taken into account.

\subsubsection{Constraints: From the Edges of the Coarse Face to the Outer Edges of the Refined Face}
We begin by adjusting the signs of sub-constraint matrices that describe the mapping from edges on the coarse
element to outer edges (specifically edges $E_4^R$ through $E_{11}^R$ in Figure \ref{fig:enumeration_of_edges}) on the refined element.
By considering the vertex order, we determine the direction of the edges and modify the corresponding entries in the constraint matrix.

\subsubsection{Constraints: From the Coarse Face to the Refined Faces}
Next, we discuss how to adapt the constraint matrix for that map to the refined faces
$F^R_0,~\dots,~F^R_3$. For an edge, there are only two possible configurations (pointing from the
left to the right or vice versa). However, in the three-dimensional case, we must consider
the $x$-direction and the $y$-direction and which direction is prioritized.
This results in $2^3=8$ possible orientations. 
Geometrically, we interpret the necessary operations as $x$-axis inversion, $y$-axis inversion, and $x$- and $y$-axis exchange.
These operations are visualized in Figure \ref{fig:face_orientation}.

\begin{figure}
  \begin{tikzpicture}[scale=0.53, transform shape, >={Stealth[length=2.5mm, width=1.5mm, inset=0.7mm]}]
    % Reference Cell
    % Draw the square
    \draw[line width=0.4mm] (0,0) -- (4,0) -- (4,4) -- (0,4) -- cycle;

    % Label the vertices
    \node at (0,0) [below left]  {\large $v_0$};
    \node at (4,0) [below right] {\large $v_1$};
    \node at (0,4) [above left]  {\large $v_2$};
    \node at (4,4) [above right] {\large $v_3$};

    % Vertices
    \fill (0,0) circle (3pt);
    \fill (0,4) circle (3pt);
    \fill (4,0) circle (3pt);
    \fill (4,4) circle (3pt);

    % Draw the arrows
    \draw[->, line width=0.0mm] (3.49, 0.20)--(3.50, 0.20);
    \draw[line width=0.2mm] (0.50, 0.25)--(3.20, 0.25);
    \draw[line width=0.2mm] (0.50, 0.15)--(3.20, 0.15);

    \draw[->, line width=0.0mm] (3.49, 3.80)--(3.50, 3.80);
    \draw[line width=0.2mm] (0.50, 3.75)--(3.20, 3.75);
    \draw[line width=0.2mm] (0.50, 3.85)--(3.20, 3.85);

    \draw[->, line width=0.2mm] (0.2, 0.50)--(0.20, 3.50);

    \draw[->, line width=0.2mm] (3.8, 0.50)--(3.80, 3.50);

    \filldraw[thin, fill=gray, fill opacity=0.25] (0,0) -- (4,0) -- (4,4) -- (0,4) -- cycle;

    % Arrow in between
    \draw[->, line width=0.3mm] (5, 2) -- (7, 2);
    \node at (6,2) [below]  {\large $x$-flip};

    % x-flip
    % Draw the square
    \draw[line width=0.4mm] (8,0) -- (12,0) -- (12,4) -- (8,4) -- cycle;

    % Label the vertices
    \node at ( 8,0) [below left]  {\large $v_1$};
    \node at (12,0) [below right] {\large $v_0$};
    \node at ( 8,4) [above left]  {\large $v_3$};
    \node at (12,4) [above right] {\large $v_2$};

    % Vertices
    \fill ( 8,0) circle (3pt);
    \fill ( 8,4) circle (3pt);
    \fill (12,0) circle (3pt);
    \fill (12,4) circle (3pt);

    % Draw the arrows
    \draw[<-, line width=0.0mm] (8.50, 0.20)--(8.51, 0.20);
    \draw[line width=0.2mm] (11.50, 0.25)--(8.70, 0.25);
    \draw[line width=0.2mm] (11.50, 0.15)--(8.70, 0.15);

    \draw[<-, line width=0.0mm] (8.50, 3.80)--(8.51, 3.80);
    \draw[line width=0.2mm] (11.50, 3.75)--(8.70, 3.75);
    \draw[line width=0.2mm] (11.50, 3.85)--(8.70, 3.85);

    \draw[->, line width=0.2mm] (8.2, 0.50)--(8.20, 3.50);

    \draw[->, line width=0.2mm] (11.8, 0.50)--(11.80, 3.50);

    \draw[dotted, thick] (10.0, -0.5) -- (10.0, 4.5);

    \filldraw[thin, fill=gray, fill opacity=0.25] (8,0) -- (12,0) -- (12,4) -- (8,4) -- cycle;

    % Arrow in between
    \draw[->, line width=0.3mm] (13, 2) -- (15, 2);
    \node at (14,2) [below]  {\large $y$-flip};

    % y-flip
    % Draw the square
    \draw[line width=0.4mm] (16,0) -- (20,0) -- (20,4) -- (16,4) -- cycle;

    % Label the vertices
    \node at (16,0) [below left]  {\large $v_3$};
    \node at (20,0) [below right] {\large $v_2$};
    \node at (16,4) [above left]  {\large $v_1$};
    \node at (20,4) [above right] {\large $v_0$};

    % Vertices
    \fill (16,0) circle (3pt);
    \fill (16,4) circle (3pt);
    \fill (20,0) circle (3pt);
    \fill (20,4) circle (3pt);

    % Draw the arrows
    \draw[->, line width=0.2mm] (19.50, 0.20)--(16.50, 0.20);
                                                                                                        
    \draw[->, line width=0.2mm] (19.50, 3.80)--(16.50, 3.80);

    \draw[->, line width=0.0mm] (16.20, 0.51)--(16.20, 0.49);
    \draw[line width=0.2mm] (16.25, 0.70)--(16.25, 3.50);
    \draw[line width=0.2mm] (16.15, 0.70)--(16.15, 3.50);

    \draw[->, line width=0.0mm] (19.8, 0.51)--(19.80, 0.49);
    \draw[line width=0.2mm] (19.85, 0.70)--(19.85, 3.50);
    \draw[line width=0.2mm] (19.75, 0.70)--(19.75, 3.50);

    \draw[dotted, thick] (15.5, 2.0) -- (20.5, 2.0);

    \filldraw[thin, fill=gray, fill opacity=0.25] (16,0) -- (20,0) -- (20,4) -- (16,4) -- cycle;

    % Arrow in between
    \draw[->, line width=0.3mm] (21, 2) -- (23, 2);
    \node at (22,2) [below]  {\large $xy$-flip};

    % xy-flip
    % Draw the square
    \draw[line width=0.4mm] (24,0) -- (28,0) -- (28,4) -- (24,4) -- cycle;

    % Label the vertices
    \node at (24,0) [below]       {\large $v_3$};
    \node at (28,0) [below right] {\large $v_1$};
    \node at (24,4) [above left]  {\large $v_2$};
    \node at (28,4) [above]       {\large $v_0$};

    % Vertices
    \fill (24,0) circle (3pt);
    \fill (24,4) circle (3pt);
    \fill (28,0) circle (3pt);
    \fill (28,4) circle (3pt);

    % Draw the arrows
    \draw[<-, line width=0.0mm] (24.50, 0.20)--(24.51, 0.20);
    \draw[line width=0.2mm] (24.70, 0.25)--(27.50, 0.25);
    \draw[line width=0.2mm] (24.70, 0.15)--(27.50, 0.15);

    \draw[<-, line width=0.0mm] (24.50, 3.80)--(24.51, 3.80);
    \draw[line width=0.2mm] (24.70, 3.75)--(27.50, 3.75);
    \draw[line width=0.2mm] (24.70, 3.85)--(27.50, 3.85);

    \draw[->, line width=0.2mm] (24.20, 3.50)--(24.2, 0.50);
                                                                                                       
    \draw[->, line width=0.2mm] (27.80, 3.50)--(27.8, 0.50);

    \draw[dotted, thick] (23.5, -0.5) -- (28.5, 4.5);

    \filldraw[thin, fill=gray, fill opacity=0.25] (24,0) -- (28,0) -- (28,4) -- (24,4) -- cycle;

\end{tikzpicture}
  \caption{
    \label{fig:face_orientation}
    Visualization of the different orientations for adapting the constraint matrix from the coarse face to the refined faces.
    We start with the reference cell and apply the $x$-axis inversion.
    Then, we apply the $y$-axis inversion to the result, followed by the $x$- and $y$-axis exchange.
  }
\end{figure}

Algorithm \ref{alg:x-inversion} demonstrates how to perform an $x$-inversion on the constraint matrix for a given cell $K$.
The $y$-inversion follows a similar approach. Additionally, Algorithm \ref{alg:xy-inversion} explains the $x$- and $y$-axis exchange.

\begin{algorithm}
  \caption{
    \label{alg:x-inversion}
    Description of the $x$-axis inversion
  }

  \tcp{Convert the double indices from the faces into one index}
  \SetKwProg{Fn}{face\_index}{:}{}
  \Fn{{$(lx,~ly$)}}{
    \Return $(lx \cdot (p - 1)) + ly$\;
  }

  \LoopOverAll ( refined face $F^R_k$ of face $F^C$ \textbf{do}) {
    \If{ The $x$-orientation of $F^R_k$ and the $x$-orientation of $F^C$ do not match } {
      \tcp{Extract the submatrix of the constraint matrix that maps the DoFs from $F^R$ onto  $F^C$.
      This corresponds to $C_{4,12+k}$ from Figure \ref{fig:constraint_matrix_structure}, where $k \in \{0,1,2,3\}$ is the number of the refined face.}
      local\_constraint\_matrix $\gets$ \textbf{get\_local\_constraint\_matrix}($F^R_k$)\;

      \tcp{Loop over the indices $i=(ix,iy)$ and $j=(jx,jy)$}
      \For{$ix = 0$ \KwTo $p-1$ } {
        \For{$iy = 0$ \KwTo $p-1$ } {
          \For{$jx = 0$ \KwTo $p-1$ } {
            \For{$jy = 0$ \KwTo $p-1$ } {
              \If{ \textbf{is\_odd}($ix + jx$) } {
                local\_constraint\_matrix(\textbf{face\_index}($ix$, $iy$), \textbf{face\_index}($jx$, $jy$)) $\gets$ -local\_constraint\_matrix(\textbf{face\_index}($ix$, $iy$), \textbf{face\_index}($jx$, $jy$))\;
              }
            }
          }
        }
      }
      \tcp{Write the modified local sub-constraint matrix into the global constraint matrix}
      \textbf{local\_constraint\_to\_global}(local\_constraint\_matrix);
    } 
  }
\end{algorithm}

\begin{algorithm}
  \caption {
    \label{alg:xy-inversion}
    Description of the $x$- and $y$-axis exchange
  }
  \tcp{Convert the double indices from the faces into one index}
  \SetKwProg{Fn}{face\_index}{:}{}
  \Fn{{$(lx,~ly$)}}{
    \Return $(lx \cdot (p - 1)) + lx$\;
  }

  \LoopOverAll ( refined face $F^R_k$ of face $F^C$ \textbf{do}) {
    \If{ The primary direction of $F^R_k$ and primary direction of $F^C$ do not match } {
      \tcp{Extract the submatrix of the constraint matrix that maps the DoFs from $F^R_k$ onto  $F^C$.}
      new\_constraint\_matrix $\gets$ \textbf{get\_local\_constraint\_matrix}($F^R_k$)\;
      old\_constraint\_matrix $\gets$ \textbf{get\_local\_constraint\_matrix}($F^R_k$)\;

      \For{$ix = 0$ \KwTo $p-1$ } {
        \For{$iy = 0$ \KwTo $p-1$ } {
          \For{$jx = 0$ \KwTo $p-1$ } {
            \For{$jy = 0$ \KwTo $p-1$ } {
              \tcp{Swap the $x$ and $y$ direction}
              new\_constraint\_matrix(\textbf{face\_index}($ix$, $iy$), \textbf{face\_index}($jx$, $jy$)) $\gets$ old\_constraint\_matrix(\textbf{face\_index}  ($ix$, $iy$), \textbf{face\_index}($jy$, $jx$))\;
            }
          }
        }
      }
      \tcp{Write the modified local sub-constraint matrix into the global constraint matrix}
      \textbf{local\_constraint\_to\_global}(new\_constraint\_matrix)
    } 
  }
\end{algorithm}

\subsubsection{Constraints: From the Edges of the Coarse Face to the Inner Edges of the Refined Face}
Next, we describe the process of adapting the constraint matrix for the inner edges, which correspond to the
edges $E_0^R, \dots, E^R_3$ from Figure \ref{fig:enumeration_of_edges}.
Finally, the most complex case is addressed last.
The constraints of the inner edges depend on all edges of the coarse parent face and the parent face itself,
as shown in Figure \ref{fig:constraint_matrix_structure}.

For the sub-constraint matrices that map from the coarse edges parallel to the refined edge, we employ
the same approach as for the outer edges. Next, we need to consider the direction of the internal edge,
which can be either in the $x$- or $y$-direction. We must apply the corresponding axis inversion,
as described above, based on the orientation of the internal edge we are currently considering.
However, we encounter an additional case for the inner edges: the sub-constraint matrix mapping
from the coarse edges orthogonal to the refined internal edge. This situation is special because,
unlike other cases, we only have the orientation of the coarse edge.
Therefore, we must test whether this orientation matches the orientation of the canonical
coarse-fine mapping or not.
We cover this by Algorithm \ref{alg:internal_edge}.

\begin{algorithm}
  \caption{
    \label{alg:internal_edge}
    Description of the inversion of the direction of the refined internal edge parallel to the $x-$axis.
  }
  \LoopOverAll ( internal edges $E^R_k$ of face $F^C$ \textbf{do}) {
    \If{ The $x$-orientation of $F^C$ differs from the canonical orientation } {
      \tcp{Extract the submatrix of the constraint matrix that maps the DoFs from $E_k^R$
           to the corresponding DoFs of $K$. Specifically, this corresponds to elements 
           $C_{k,2}$ and $C_{k,3}$ in Figure \ref{fig:constraint_matrix_structure}. 
           Notice: We need to perform this operation twice, once for each submatrix.}
      local\_constraint\_matrix $\gets$ \textbf{get\_local\_constraint\_matrix}($E^R_k$)\;

      \For{$i,j$ in local\_constraint\_matrix } {
        \If{ \textbf{is\_odd}(i + j) } {
          local\_constraint\_matrix(i, j) $\gets$ -local\_constraint\_matrix(i, j)\;
        }
      }
      \tcp{Write the modified local sub-constraint matrix into the global constraint matrix}
      \textbf{local\_constraint\_to\_global}(local\_constraint\_matrix);
    } 
  }
\end{algorithm}

%+-------------------------------------------------------------------+
%|             Section: Model Problem and Numerical Tests            |
%+-------------------------------------------------------------------+
\section{Model Problem and Numerical Tests}
\label{sec:numerical_tests}
In this section, we introduce the time-harmonic Maxwell's equations as a model problem.
We present two numerical examples demonstrating
our implementation of hanging nodes for Nédélec elements, especially for non-orientable locally refined meshes.

Let $\Omega \subset \mathbb{R}^3$
be a bounded Lipschitz domain with a sufficiently smooth boundary
$\Gamma = \Gamma^{\text{inc}} \cup \Gamma^{\infty}$, where on $\Gamma^{\infty}$ an 
absorbing boundary condition is given
and on $\Gamma^{\text{inc}}$, a boundary condition for some given incident electric field is given.
Find the electric field ${\vec u} \in \mathbf{H}_\text{curl}(\Omega)$ such that for all $\vec \varphi \in \mathbf{H}_\text{curl}(\Omega)$ it holds
\begin{equation}
  \label{eq:THM_weak}
  \begin{split}
    \int_\Omega \left( \mu^{-1} \operatorname{curl} \left( {\vec u} \right) \cdot \operatorname{curl} \left( \vec \varphi \right)
    - \varepsilon \omega^2 {\vec u} \cdot\vec \varphi \right) ~ \mathsf{d} x
    + i \kappa \omega \int_{\Gamma} 
    (\vec{n}\times (\vec{u} \times \vec{n}))\cdot (\vec{n}\times (\vec{\varphi}\times\vec{n} )) ~ \mathsf{d} s \\
    = \int_{\Gamma^\text{inc}} 
    (\vec{n}\times (\vec{u}^\text{inc} \times \vec{n})) \cdot (\vec{n}\times (\vec{\varphi}\times \vec{n} )) ~ \mathsf{d} s 
  \end{split}
\end{equation}
with the outer normal vector $\vec{n}$.
Here, ${\vec u}^\text{inc}$ with $ \vec{n}\times {\vec u}^\text{inc}\in L^2(\Gamma^{\text{inc}},\mathbb{C}^d)$
is some given incident electric field, 
$\mu \in \mathbb{R}^+$ is the relative magnetic permeability, $\kappa = \sqrt{\varepsilon}$, 
$\varepsilon \in \mathbb{R}^+$ is the relative permittivity, $\omega = \frac{2 \pi}{\lambda}$ is the wavenumber, and $\lambda \in \mathbb{R}^{+}$ is the wavelength. 
System~\eqref{eq:THM_weak} is called time-harmonic because the time dependence can be expressed by $e^{i \omega \tau}$,
where $\tau > 0$ denotes the time. For the derivation of the time-harmonic Maxwell's equations,
we refer the reader to \cite{Bk:Monk:2003}.

We briefly comment on the numerical solution of the resulting linear systems, which is rather challenging as it is ill-posed. Consequently, specialized methods have to be employed. 
A well-known approach to address the time-harmonic Maxwell's equation is based on combining direct solvers and domain decomposition methods \cite{Art:Gander:OptimizedSchwarz:12,Art:Kinnewig:DD26:21}. 
Here, the basic idea is to divide the problem into small enough sub-problems 
so that a direct solver can handle each sub-problem.
Another approach is to find suitable preconditioners for iterative solvers, for example, with the help
of $\mathcal{H}$-matrices \cite{Art:FauParMel:23}. 
As the computation of such preconditioners is quite challenging, these methods can be
combined with a domain decomposition method \cite{Art:Parvizi:23}.

\subsection{Qualitative and Quantitative Computational Analysis on a Simple Waveguide}
\label{subsec:quantitaive_analysis}
In this first numerical example, we investigate qualitatively, in terms of the `picture norm', as well quantitatively, in terms of a small convergence analysis on a sequence of locally refined meshes, our newly proposed algorithms, and implementation.
We consider a simplified model of glass fiber, 
which is modeled by the domain $\Omega = (0, 4) \times (0, 4) \times (0, 1.5) \; \mu m$ with a cylindrical structure in the center. 
The center is made from SiO$_2$ with a refractive index of
$n_{\text{SiO}_2}=2.0257$ ($\mu_{\text{SiO}_2} = 1.0000$, $\varepsilon_{\text{SiO}_2} = n_{\text{SiO}_2}^2$) 
surrounded by air $n_\text{air}=1.0000$ ($\mu_{\text{air}} = 1.0000$, $\varepsilon_{\text{air}} = 1.0000$), 
an incident wave with a wavelength of $\lambda=375\;nm$. 
The geometry is shown in Figure~\ref{fig:geometry_quantitative_analysis} (left).
The incoming electric field is represented by $u_{\text{inc}} = \operatorname{exp}\
\left(\frac{-20}{\mu m^2} (x^2+y^2)\right) \vec{e}_x$ with unit vector $\vec{e}_x$ in $x$-direction.
Furthermore, $\Gamma_{\text{inc}} =  (0, 4) \times (0, 4) \times \{ 0 \} \; \mu m$ 
denotes the boundary with the incident boundary condition, while all other boundaries $\Gamma_\infty$ are characterized by absorbing conditions, namely homogeneous Robin conditions.

We evaluate the following three goal functionals: 
the point value $J_P(u) = u_{0} \left( P \right)$, where $P = (2.2\;\mu m, 2.2\;\mu m, 0.2\;\mu m)$, 
the face integral $J_F(u) = \left\lVert (\vec u - \vec u_\text{ref}) \times \vec n \right\rVert_{L^2({\Gamma}_\text{out})}$ 
where ${\Gamma}_\text{out} =  (0, 4) \times (0, 4) \times \{ 1.5 \} \; \mu m$ 
and the domain integral $ J_D(u) = \left\lVert (\vec u - \vec u_\text{ref}) \right\rVert_{L^2(\Omega)}$.
On the finest level with $2\, 080\, 944$ DoFs, the numerical solution is used as the numerical reference value.
The results are presented in Table~\ref{tab:quantitative_tests}.
In this test, we employ the polynomial degree of the underlying base functions high enough, namely $p=3$, so that all features of the base functions are tested.

% TODO: Still in the TOMS style
\begin{table}[h]
  \begin{center}
    \begin{tabular}{ l l l l l } \toprule
    Level $l$ & DoFs   & $|J_P(u_l) - J_P(u_{ref})|$ & $|J_F(u_l)|$ & $|J_D(u_l)|$ \\ \midrule
    1         & 29436  & 0.052260                    & 0.00560482   & 0.000530526 \\
    2         & 146520 & 0.010761                    & 0.00316504   & 0.000254373 \\
    3         & 681432 & 0.000079                    & 0.00166421   & 0.000152532 \\ \bottomrule
    \end{tabular}
    \caption{Section \ref{subsec:quantitaive_analysis}. Results from evaluating the goal functionals on different levels.}
    \label{tab:quantitative_tests}
  \end{center}
\end{table}

\begin{figure}
  \begin{center}
    \begin{minipage}{0.45\textwidth}
      \begin{tikzpicture}[x={(1cm,0cm)}, y={(0cm,1cm)}, z={(0.5cm,0.5cm)}]
    % Cube
    \draw[white] (4.0, 4.0, 0.0) -- (4.0, 4.0, 3.2);

    \filldraw[fill=black!50, thick] (4.0, 0.0, 0.0) -- (4.0, 0.0, 1.5) -- (4.0, 4.0, 1.5) -- (4.0, 4.0, 0.0) -- cycle;
    \filldraw[fill=black!30, thick] (0.0, 0.0 ,0.0) -- (4.0, 0.0, 0.0) -- (4.0, 4.0, 0.0) -- (0.0, 4.0, 0.0) -- cycle;
    \filldraw[fill=black!50, thick] (0.0, 4.0, 0.0) -- (0.0, 4.0, 1.5) -- (4.0, 4.0, 1.5) -- (4.0, 4.0, 0.0) -- cycle;

    \draw[dotted] (0.0, 0.0, 0.0) -- (0.0, -1.0, 0.0);
    \draw[dotted] (4.0, 0.0, 0.0) -- (4.0, -1.0, 0.0);

    \draw[dotted] (-1.0, 0.0, 0.0) -- (0.0, 0.0, 0.0);
    \draw[dotted] (-1.0, 4.0, 0.0) -- (0.0, 4.0, 0.0);

    \draw[dotted] (4.0, 0.0, 1.5) -- (4.0, -1.0, 1.5);

    \draw[<->] (-0.5, 0.0, 0.0) -- (-0.5, 4.0, 0.0) node[midway, left] {$4\mu m$};
    \draw[<->] (0.0, -0.5, 0.0) -- (4.0, -0.5, 0.0) node[midway, below] {$4\mu m$};
    \draw[<->] (4.0, -0.5, 0.0) -- (4.0, -0.5, 1.5) node[below, yshift=-10] {$1.5\mu m$};

    \filldraw[fill=black!10] (2,2) circle (1.4);
    \draw[<->] (0.6, 2.0, 0.0) -- (3.4, 2.0, 0.0) node[midway, below] {$2.8\mu m$};

\end{tikzpicture}
    \end{minipage}
    \hfill
    \begin{minipage}{0.45\textwidth}
      \begin{tikzpicture}[x={(1cm,0cm)}, y={(0cm,1cm)}, z={(0.25cm,0.25cm)}]
    % Cube
    \filldraw[fill=black!50, thick] (4.0, 0.0, 0.0) -- (4.0, 0.0, 6.4) -- (4.0, 4.0, 6.4) -- (4.0, 4.0, 0.0) -- cycle;
    \filldraw[fill=black!30, thick] (0.0, 0.0 ,0.0) -- (4.0, 0.0, 0.0) -- (4.0, 4.0, 0.0) -- (0.0, 4.0, 0.0) -- cycle;
    \filldraw[fill=black!50, thick] (0.0, 4.0, 0.0) -- (0.0, 4.0, 6.4) -- (4.0, 4.0, 6.4) -- (4.0, 4.0, 0.0) -- cycle;

    \draw[dotted] (0.0, 0.0, 0.0) -- (0.0, -1.0, 0.0);
    \draw[dotted] (4.0, 0.0, 0.0) -- (4.0, -1.0, 0.0);

    \draw[dotted] (-1.0, 0.0, 0.0) -- (0.0, 0.0, 0.0);
    \draw[dotted] (-1.0, 4.0, 0.0) -- (0.0, 4.0, 0.0);

    \draw[dotted] (4.0, 0.0, 6.4) -- (4.0, -1.0, 6.4);

    \draw[<->] (-0.5, 0.0, 0.0) -- (-0.5, 4.0, 0.0) node[midway, left] {$16\mu m$};
    \draw[<->] (0.0, -0.5, 0.0) -- (4.0, -0.5, 0.0) node[midway, below] {$16\mu m$};
    \draw[<->] (4.0, -0.5, 0.0) -- (4.0, -0.5, 6.4) node[below, yshift=-14] {$25.6\mu m$};

    \filldraw[fill=black!10] (2.0, 3.2, 0.0) circle (0.28);
    \filldraw[fill=black!10] (2.0, 0.8, 0.0) circle (0.28);

    \filldraw[fill=black!10] (0.8, 2.6, 0.0) circle (0.28);
    \filldraw[fill=black!10] (0.8, 1.4, 0.0) circle (0.28);

    \filldraw[fill=black!10] (3.2, 2.6, 0.0) circle (0.28);
    \filldraw[fill=black!10] (3.2, 1.4, 0.0) circle (0.28);

    %\draw[<->] (2, 0.8, 0.0) -- (3.2, 1.4, 0.0) node[midway, below, xshift=5] {$3\mu m$};
    \draw[<->] (0.8, 1.4, 0.0) -- (0.8, 2.6, 0.0) node[midway, right] {$3\mu m$};

    \draw[densely dotted] (2.92, 2.6, 0.0) -- (2.92, 3.00, 0.0);
    \draw[densely dotted] (3.48, 2.6, 0.0) -- (3.48, 3.00, 0.0);
    \draw[<->] (2.92, 2.95, 0.0) -- (3.48, 2.95, 0.0) node[midway, above] {$1.4\mu m$};

\end{tikzpicture}
    \end{minipage}
    \caption{
      Left: Section \ref{subsec:quantitaive_analysis}. Geometry and dimensions of the simplified waveguide.
      The core, marked in light gray, is made of SiO$_2$ which is surrounded by air (dark grey).
      Right: Section \ref{sec_waveguide}.  Geometry and dimensions of the waveguide. In the 
      so-called modifications (light gray), the refractive index is higher than in the surroundings.
    }
    \label{fig:geometry_quantitative_analysis}
  \end{center}
\end{figure}

\begin{figure}
  \begin{tabular}{ccccc} \toprule
    &FE\_Nédélec & FE\_NédélecSZ & FE\_NédélecSZ extended & \\ \midrule
    \begin{sideways}~~~~~~~~~~~Level 1\end{sideways} &
    \includegraphics[width=0.25\textwidth]{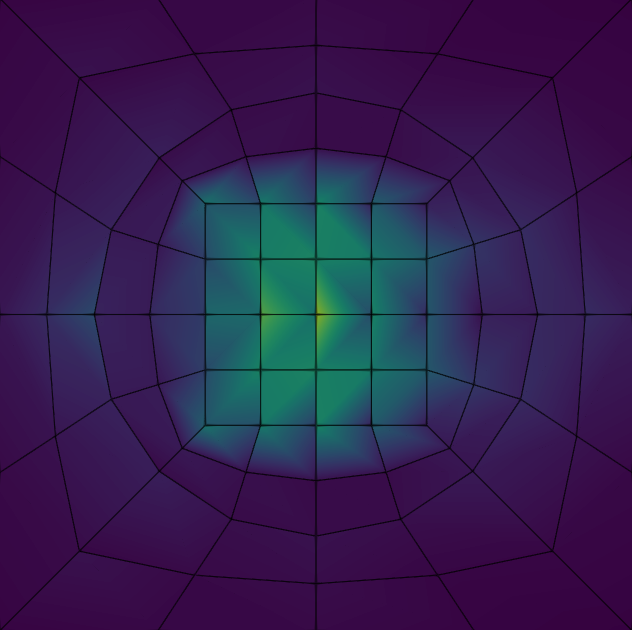} &
    \includegraphics[width=0.25\textwidth]{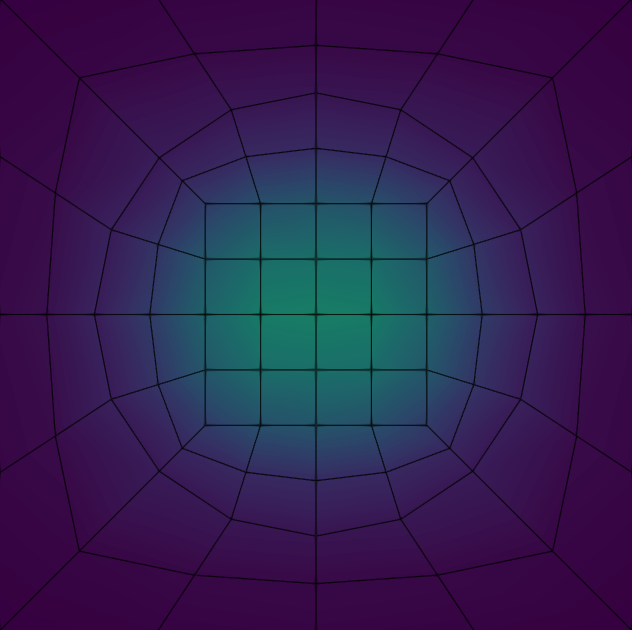} &
    \includegraphics[width=0.25\textwidth]{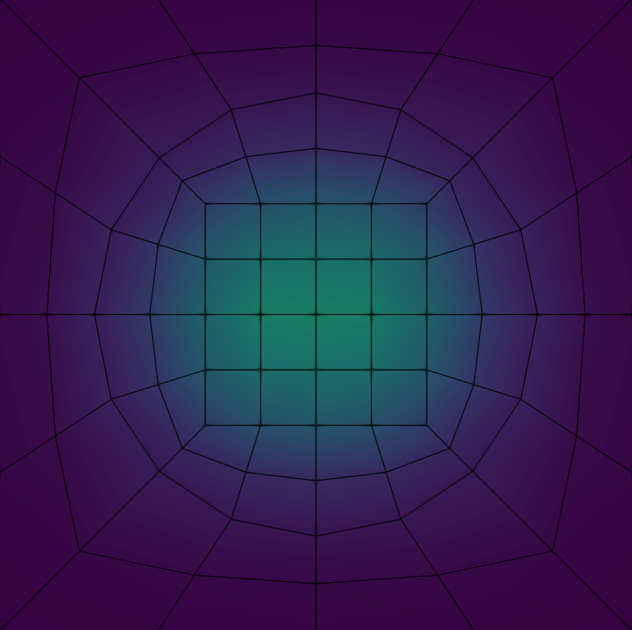} &
    \includegraphics[width=0.08\textwidth]{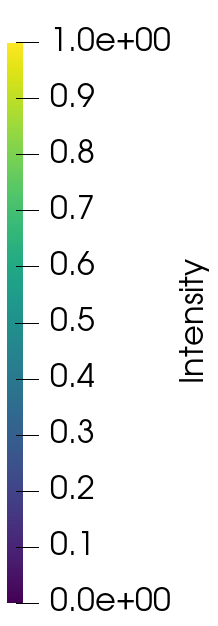} \\
    \begin{sideways}~~~~~~~~~~~Level 3\end{sideways} &
    \includegraphics[width=0.25\textwidth]{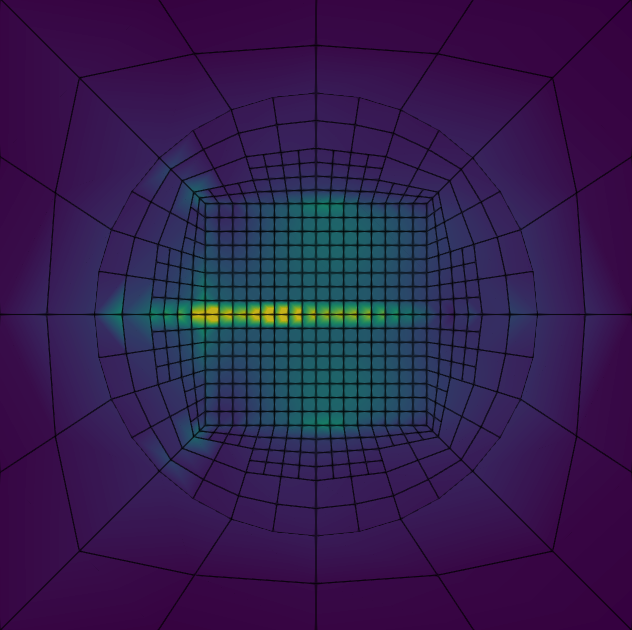} &
    \includegraphics[width=0.25\textwidth]{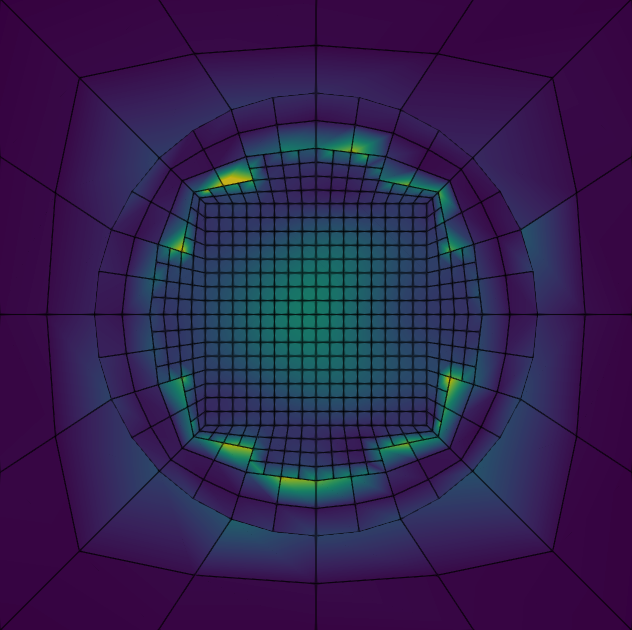} &
    \includegraphics[width=0.25\textwidth]{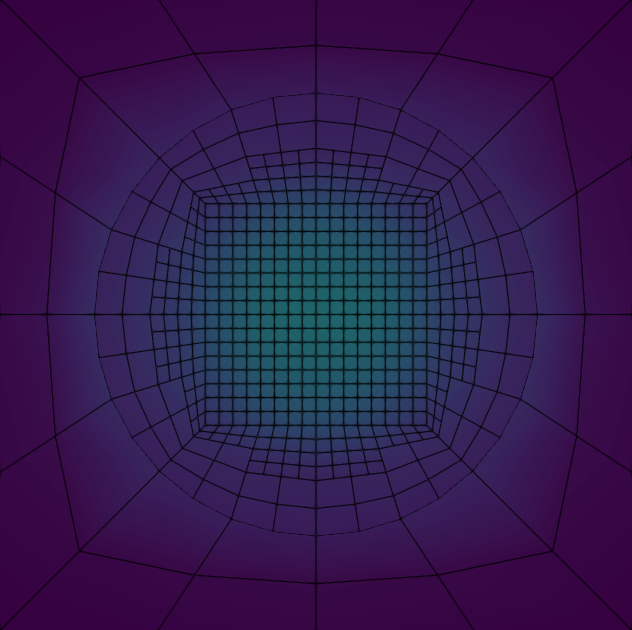} &
    \includegraphics[width=0.08\textwidth]{Images/Nedelec-Color-Key} \\ \bottomrule
  \end{tabular}
  \caption{%
    Section \ref{subsec:quantitaive_analysis}.
    Comparison in the `picture norm' of the different implementations from the N\'ed\'elec elements on the example
    of the intensity plot from the fiber
    for two refinement levels. Uniform refinement was applied to the first row. Local mesh refinement causing hanging edges is shown in the second row. In the columns, we have
    from left to right FE\_N\'ed\'elecSZ, FE\_N\'ed\'elecSZ, and our newly proposed extension of FE\_N\'ed\'elecSZ.
    We clearly observe the wrong implementations of FE\_N\'ed\'elec (left) and FE\_N\'ed\'eleccSZ (middle).
  }
  \label{fig:comparision}
\end{figure}

In Figure~\ref{fig:comparision}, we compare against the existing implementation of the N\'ed\'elec elements in \dealii, where the errors resulting from the sign conflict are visible. 
The plots in the first column are computed using the FE\_N\'ed\'elec class, which does not support non-oriented meshes. 
Therefore, the resulting intensity distribution differs from the correct solution.

The results computed with the existing implementation of the FE\_N\'ed\'elecSZ class are presented in the
second column.
Here, the solution on the uniform refined grid is correct, but on the isotropic refined grid, the solution
differs from the correct solution.
The results from our proposed extension of the FE\_N\'ed\'elecSZ class are shown in the third column.
Specifically, the numerical solution on both grids (locally refined and uniformly refined) is correct.

% TODO: Still in TOMS Style

\subsection{Laser-Written Waveguide}
\label{sec_waveguide}
As a second example, a practical application in optics simulations is considered. 
To guide optical waves, we need a difference in the refractive index. 
This can be achieved by causing stress and compression in the material. 
These changes (modifications) can be introduced by hitting the material with a femtosecond laser pulse, 
creating a quickly expanding plasma and introducing stress and compression. 
Here, the modifications form a hexagonal pattern, making the material denser in its center, 
leading to a contrast in the refractive index. 
For a more detailed description of the geometry and the process of creating such waveguides,
we refer the reader to \cite{Art:Perevoznik:21}.

We consider the domain $\Omega = (0.0, 16.0) \times (0.0, 16.0) \times (0.0, 25.6) \; \mu m$ shown in Figure~\ref{fig:geometry_quantitative_analysis} (right). Here, we have
the incident boundary $\Gamma_{\text{inc}} = (0.0, 16.0) \times (0.0, 16.0) \times \{0\} \; \mu m$ and the 
incident electric field 
\[
  u_{\text{inc}} = \operatorname{exp}\left(\frac{-57}{\mu m^2} \left((x-0.5\;\mu m)^2+(y-0.5\;\mu m)^2\right)\right) \vec{e}_y.
\]
All other boundaries $\Gamma_\infty$ are absorbing boundaries, i.e., homogeneous Robin boundaries.
Concerning the material properties, we assume the carrier material to have a refractive index of 
$n_\text{cladding}=1.4899$ ($\mu_\text{cladding} = 1.0000$, $\varepsilon_\text{cladding} = n_\text{cladding}^2$),
the compressed center to have a refractive index of 
$n_\text{center}=1.4906$ ($\mu_\text{center} = 1.0000$, $\varepsilon_\text{center} = n_\text{center}^2$)
and the modifications to be filled with air (see section \ref{subsec:quantitaive_analysis}).
The incident laser light has a wavelength of $\lambda=660\;nm$.

\begin{figure}[h]
  \begin{center}
    \includegraphics[width=0.76\textwidth]{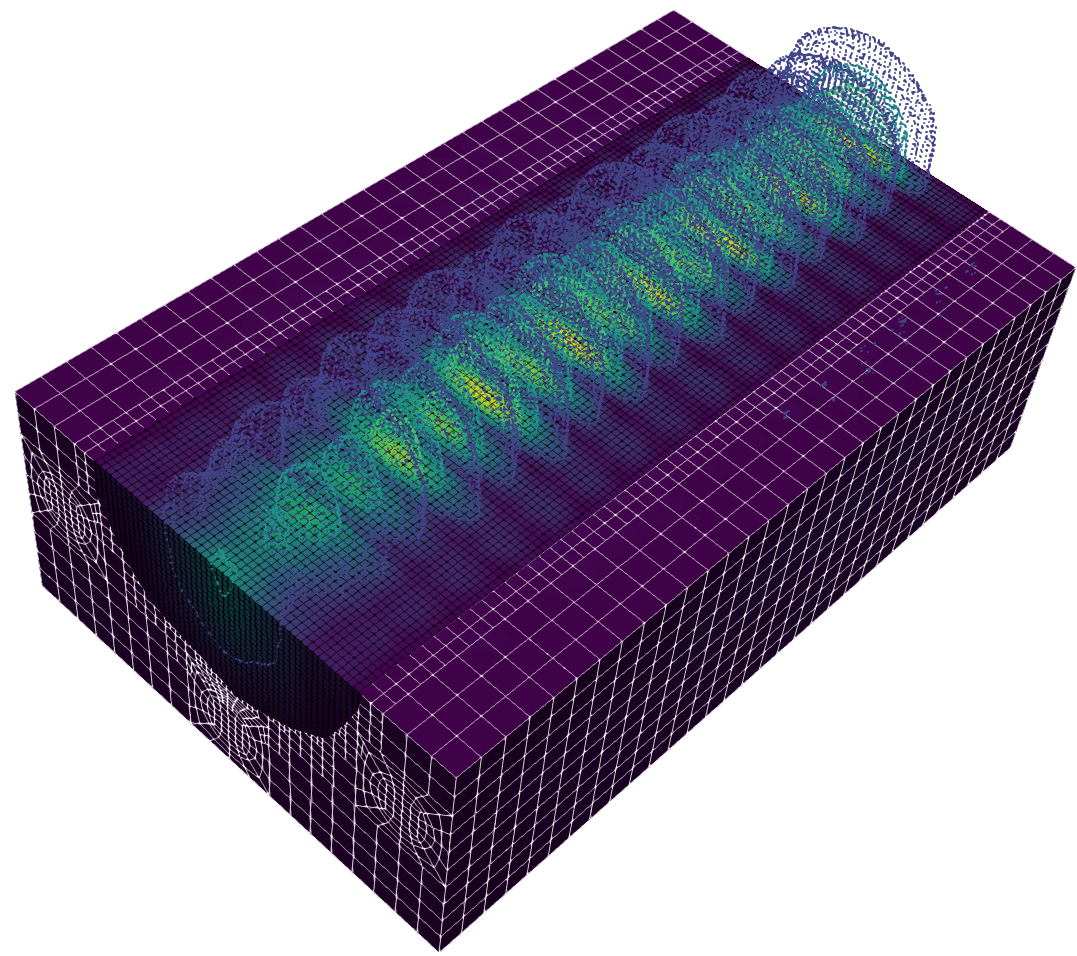}
  \end{center}
  \caption{%
    Section \ref{sec_waveguide}: Cross-section through the waveguide at the plane $(0, 16) \times \{8\} \times (0,25.6) \; \mu m$,
    where the electric field's intensity distribution inside the waveguide is visualized.
    The black edges represent the edges of the finest level, while the white lines show
    the edges of the coarser levels.
    }
  \label{fig_example_2}
\end{figure}

As the geometry is rather complex, the numerical efforts of such a three-dimensional configuration in terms 
of computational cost require the domain decomposition implemented in \cite{Art:Kinnewig:DD26:21}, 
where $48$ subdomains were employed and local mesh refinement as shown in Figure \ref{fig_example_2}.

A discussion and interpretation of this example is as follows. As previously mentioned, we deal with hanging edges and faces due to local mesh refinement. Without our newly proposed extensions of FE\_N\'ed\'elecSZ,
such computations on complex geometries from practical applications would not have been possible and demonstrate the capabilities of both the algorithmic advancements in this work as well as our open-source codes\footnote{\url{https://zenodo.org/records/10913219}} in \dealii.

%+-------------------------------------------------------------------+
%|               Section: Conclusion                                 |
%+-------------------------------------------------------------------+
\section{Conclusion}
In this work, we have addressed the sign conflict problem  in three spatial dimensions of the N\'ed\'edec elements that 
appear in scenarios where hanging nodes arise on locally refined meshes.
We provided a comprehensive derivation
in terms of algorithmic designs for resolving this sign conflict. These concepts can be applied to any software package that supports N\'ed\'elec elements and locally refined meshes on
quadrilaterals or hexahedra with hanging nodes.
Our choice of \dealii as a programming platform has proven to be highly accessible and user-friendly.
The new implementation was demonstrated for 
two
numerical experiments that include qualitative comparisons in
three spatial dimensions as well as some
computational convergence studies.
In the second numerical example, we presented a 
practical example from optics simulations showing a laser-written waveguide.
Not only does this example validate our implementation of hanging nodes for Nédélec elements on non-orientable grids, but it
also demonstrates its practical application in optics simulations on complex geometries where local mesh refinement is indispensable.

%+-------------------------------------------------------------------+
%|               Section: Acknowledgements                           |
%+-------------------------------------------------------------------+
\section*{Acknowledgments}
  This work is funded by the Deutsche Forschungsgemeinschaft (DFG) under Germany’s Excellence
  Strategy within the Cluster of Excellence PhoenixD (EXC 2122, Project ID 390833453).
  Furthermore, we would like to thank Tim Haubold and Philipp König for many fruitful
  discussions, Martin Kronbichler for helpful feedback on the actual implementation,
  and Clemens Pechstein for discussions on debugging the implementation.

%+-------------------------------------------------------------------+
%|               Section: References                                 |
%+-------------------------------------------------------------------+
\bibliography{lit}

\end{document}